\documentclass{article}

\usepackage{times}
\usepackage{graphicx} 
\usepackage{subfigure}

\usepackage{natbib}

\usepackage{algorithm}
\usepackage{algorithmic}

\usepackage{hyperref}

\usepackage[accepted]{icml2015}

\icmltitlerunning{An Asynchronous Distributed Proximal Gradient Method for Composite Convex Optimization}

\usepackage{amsmath, amsfonts, amssymb, amsthm,dsfont}
\usepackage{booktabs,multirow,enumerate}

\DeclareMathOperator*{\sgn}{sgn}
\newcommand{\mb}[1]{\mathbf{#1}}
\def\fal{\textbf{DFAL}}
\def\afal{\textbf{AFAL}}
\def\apg{\textbf{MS-APG}}
\def\bcd{\textbf{RBCD}}
\def\abcd{\textbf{ARBCD}}
\def\admm{\textbf{ADMM}}
\def\sadmm{\textbf{SADMM}}
\def\fprod#1{\left\langle#1\right\rangle}
\def\prox#1{\mathbf{prox}_{#1}}

\def\tk{\theta^{(k)}}
\def\lk{\lambda^{(k)}}
\def\ak{\alpha^{(k)}}
\def\xik{\xi^{(k)}}
\def\li{\lambda^{(1)}}
\def\ai{\alpha^{(1)}}
\def\xii{\xi^{(1)}}
\def\pk{P^{(k)}}
\def\fk{f^{(k)}}
\def\xk{\bx^{(k)}}
\def\gk{g^{(k)}}
\def\xkopt{\bx_*^{(k)}}
\def\Lk{L^{(k)}}
\def\Nk{N^{(k)}}
\def\T{\mathsf{T}}

\newtheorem*{mthm}{Main Result}
\newtheorem{defn}{Definition}
\newtheorem{assump}{Assumption}
\newtheorem{cor}{Corollary}

\newtheorem{thm}{Theorem}
\newtheorem{lemma}{Lemma}

\def\grad{\nabla}

\def\bh{\mathbf{h}}

\def\bq{\mathbf{q}}

\def\bu{\mathbf{u}}

\def\bx{\mathbf{x}}  
\def\by{\mathbf{y}}
\def\bz{\mathbf{z}}

\def\cC{\mathcal{C}}
\def\cD{\mathcal{D}}
\def\cE{\mathcal{E}}

\def\cG{\mathcal{G}}

\def\cI{\mathcal{I}}

\def\cK{\mathcal{K}}

\def\cN{\mathcal{N}}
\def\cO{\mathcal{O}}
\def\cP{\mathcal{P}}

\def\cR{\mathcal{R}}

\def\cY{\mathcal{Y}}
\def\cZ{\mathcal{Z}}

\def\eq{\[}
\def\en{\]}

\def\smskip{\smallskip}

\def\texitem#1{\par\smskip\noindent\hangindent 25pt
               \hbox to 25pt {\hss #1 ~}\ignorespaces}

\def\norm#1{\|#1\|}

\newcommand{\BEAS}{\begin{eqnarray*}}
\newcommand{\EEAS}{\end{eqnarray*}}
\newcommand{\BEA}{\begin{eqnarray}}
\newcommand{\EEA}{\end{eqnarray}}
\newcommand{\BEQ}{\begin{eqnarray}}
\newcommand{\EEQ}{\end{eqnarray}}
\newcommand{\BIT}{\begin{itemize}}
\newcommand{\EIT}{\end{itemize}}
\newcommand{\BNUM}{\begin{enumerate}}
\newcommand{\ENUM}{\end{enumerate}}

\newcommand{\BA}{\begin{array}}
\newcommand{\EA}{\end{array}}

\newcommand{\reals}{\mathbb{R}}
\newcommand{\integers}{\mathbb{Z}}

\newcommand{\Rank}{\mathop{\bf rank}}

\newcommand{\argmin}{\mathop{\rm argmin}}

\newif\ifpagenumbering
\pagenumberingtrue

\pagenumberingfalse

\newsavebox{\theorembox}
\newsavebox{\lemmabox}
\newsavebox{\assbox}
\savebox{\theorembox}{\noindent\bf Theorem}
\savebox{\lemmabox}{\noindent\bf Lemma}
\savebox{\assbox}{\noindent\bf Assumption}

\begin{document}

\twocolumn[
\icmltitle{An Asynchronous Distributed Proximal Gradient Method\\ for Composite Convex Optimization}

\icmlauthor{N. S. Aybat, and Z. Wang}{nsa10@psu.edu,\ zxw121@psu.edu}
\icmladdress{Penn State University,
            University Park, PA 16802 USA}
\icmlauthor{G. Iyengar}{garud@ieor.columbia.edu}
\icmladdress{Columbia University,
            New York, NY 10027 USA}

\icmlkeywords{asynchronous, distributed optimization, composite convex function, augmented Lagrangian, first-order method, spectral gap}
\vskip 0.1in
]
\begin{abstract}
\vspace{-0.25cm}
We propose a distributed first-order augmented
Lagrangian~(\fal) algorithm to minimize the sum of
composite convex functions, where each term in the sum is a private cost function belonging to a node,
and only nodes connected by an edge can directly communicate with each other. This optimization model
abstracts
a number of applications in distributed sensing and machine learning. 
We show that any limit point of \fal~iterates is optimal; and for any
$\epsilon>0$, an $\epsilon$-optimal and $\epsilon$-feasible
solution  can be computed within $\cO(\log(\epsilon^{-1}))$ \fal~iterations, which require $\cO(\frac{\psi_{\max}^{1.5}}{d_{\min}}~\epsilon^{-1})$ proximal gradient computations and
communications per node in total, where $\psi_{\max}$ denotes the largest eigenvalue
of the graph Laplacian, and $d_{\min}$ is the minimum degree of the graph.
We also propose an asynchronous version of \fal~by incorporating
randomized block coordinate descent methods; and demonstrate the
efficiency of \fal~on large scale sparse-group LASSO problems.\vspace{-0.25cm}
\end{abstract}
\section{Introduction}
\label{sec:intro}
Let $\cG=(\cN,\cE)$ denote a \emph{connected} undirected graph of $N$
computing nodes where nodes $i$ and $j$ can communicate information
only if $(i,j) \in \cE$.  Each node $i\in\cN:=\{1,\ldots,N\}$ has a
\emph{private} (local) cost function  
\vspace*{-0.25cm}
\begin{equation}
\label{eq:F_i}
F_i(x):= \rho_i(x) + \gamma_i(x), \vspace*{-0.2cm}
\end{equation}
where $\rho_i: \mathbb{R}^n \rightarrow \mathbb{R}$ is a possibly
non-smooth convex function, and $\gamma_i: \mathbb{R}^n \rightarrow
\mathbb{R}$ is a
smooth convex function. 
We assume that the proximal map \vspace*{-0.25cm}
\begin{equation}
  \label{eq:prox}
  \prox{\rho_i}(x):=\argmin_{y\in\reals^n}\left\{\rho_i(y)+\tfrac{1}{2}\norm{y-x}_2^2\right\} \vspace*{-0.25cm}
\end{equation}
is \emph{efficiently} computable for $i\in\cN$.

We propose a distributed augmented Lagrangian
algorithm for efficiently computing a solution for the
convex problem:
\vspace*{-0.35cm}
\begin{equation}\label{eq:centr}
F^*:=\min_{x \in \mathbb{R}^n}  F(x):=\sum_{i=1}^N F_i(x). \vspace*{-0.15cm}
\end{equation}
Clearly, \eqref{eq:centr} can be solved in a ``centralized" fashion by
communicating all the private
functions $F_i$ to a central node, and solving the overall problem at this
node. However, such an approach can be very expensive both from communication
and computation perspectives. Suppose $(A_i,b_i) \in \reals^{m \times (n+1)}$ and $F_i(x) = \norm{A_i x - b_i}_2^2 +\lambda \norm{x}_1$ for $i\in\cN$ such that $m\ll n$ and $N\gg 1$. Hence, \eqref{eq:centr} is a very large scale LASSO problem \emph{distributed} data. 
To solve~\eqref{eq:centr} in a centralized fashion the data
$\{(A_i,b_i): i \in
\cN\}$ needs to be communicated to the central node. This can be
prohibitively expensive, and may also violate privacy constraints. 
Furthermore, it requires that the central node have large enough memory to
be able to
accommodate all the data. On the other hand, at the expense of slower
convergence, one can completely do away with a central node, and seek for
\emph{consensus} among all the nodes on an optimal decision
using ``local" decisions communicated by the neighboring nodes. In addition, 
for certain cases, computing partial gradients \emph{locally} 
in an \emph{asynchronous} 
manner can be even more computationally efficient
when compared to computing the entire gradient at a central node. With these
considerations in mind, we
propose decentralized algorithms that can compute
solutions to~\eqref{eq:centr} using only local computations; thereby,
circumventing all privacy, communication and memory issues.
To facilitate the design of decentralized algorithms, we take advantage of
the fact that graph $\cG$ is connected, and
reformulate~\eqref{eq:centr} as \vspace*{-0.25cm}
\begin{equation}
  \label{eq:dist_opt}
  \min_{x_i \in \mathbb{R}^n,~i\in\cN}\Big\{\sum_{i=1}^N F_i(x_i):\
  x_i=x_j,\ \forall~(i,j)\in\cE\Big\}. \vspace*{-0.2cm}
\end{equation}
Optimization problems of form~\eqref{eq:dist_opt} model a variety of very
important
applications, e.g., distributed linear
regression~\cite{mateos2010distributed}, distributed
control~\cite{necoara2008application}, machine
learning~\cite{mcdonald2010distributed}, and estimation using sensor
networks~\cite{lesser2003distributed}.

We call a solution $\mathbf{\bar{x}} = (\bar{x}_i)_{i\in \cN}$ $\epsilon$-feasible
if the consensus violation
$\max_{(i,j)\in\cE}\big\{\norm{\bar{x}_i-\bar{x}_j}_2\big\} \leq \epsilon$
and $\epsilon$-optimal if $\big|\sum_{i\in\cN}F_i(\bar{x}_i)-F^*\big| \leq
\epsilon$.
In this work, we propose a \emph{distributed} first-order augmented
Lagrangian~(\fal) algorithm, establish the following main result for the synchronous case in Section~\ref{sec:laplacian}, and extend it to an \emph{asynchronous} setting in Section~\ref{sec:fal_async}.
\begin{table*}[!t]
  \centering
  {\scriptsize
    \begin{tabular}{|l|l|l|l|l|l|l|l|}
      \hline
      \multirow{2}{*}{\textbf{Reference}} & \multirow{2}{*}{\textbf{assumption on} $F_i$} & \multirow{2}{*}{operation / iter.} &iter \# for & iter \# for & comm. steps& \multirow{2}{*}{asnych.}& Can handle\\
      & & & $\epsilon$-feas. & $\epsilon$-opt. & $\epsilon$-opt. & &constraints?\\ \hline
      \citet{Duchi12} & convex, Lipschitz cont. & subgrad., projection& unknown & $\mathcal{O}(1/\epsilon^2)$ &$\cO(1/\epsilon^2)$ & no & no\\
      \hline
      \citet{nedic2009distributed} & convex & subgrad. &$\mathcal{O}(1)$ & $\mathcal{O}(1/\epsilon^2)$ & $\mathcal{O}(1/\epsilon^2)$& no & no\\
      \hline
      \citet{wei2012_1} & strictly convex & $\prox{F_i}$& unknown & $\cO(1/\epsilon)$&$\cO(1/\epsilon)$ & no & no\\ \hline
      \citet{makhdoumi2014broadcast} & convex & $\prox{F_i}$ & $\cO(1/\epsilon)$& $\cO(1/\epsilon)$& $\cO(1/\epsilon)$& no & no\\ \hline
      \citet{wei20131} & convex & $\prox{F_i}$ & $\cO(1/\epsilon)$ & $\cO(1/\epsilon)$ & $\cO(1/\epsilon)$ & yes & no\\ \hline
      \multirow{2}{*}{\citet{jakovetic2011fast}} & smooth convex& \multirow{2}{*}{$\grad{F_i}$} &\multirow{2}{*}{$\mathcal{O}(1/\sqrt{\epsilon})$} &  \multirow{2}{*}{$\mathcal{O}(1/\sqrt{\epsilon})$} &\multirow{2}{*}{$\mathcal{O}(1/\sqrt{\epsilon})$} & \multirow{2}{*}{no} & \multirow{2}{*}{no}\\
      & bounded $\grad F_i$ & & & & & &\\ \hline
      \multirow{2}{*}{\citet{chen2012fast}} & composite convex $F_i=\rho+\gamma_i$& \multirow{2}{*}{$\prox{\rho}$, $\grad \gamma_i$}& \multirow{2}{*}{$\cO(1/\sqrt{\epsilon})$} & \multirow{2}{*}{$\cO(1/\sqrt{\epsilon})$} & \multirow{2}{*}{$\cO(1/\epsilon)$}& \multirow{2}{*}{no} & \multirow{2}{*}{no}\\
      & bounded $\grad\gamma_i$ & & & & & &\\ \hline
       Our work & composite convex $F_i=\rho_i+\gamma_i$& $\prox{\rho_i}$, $\grad \gamma_i$ &$\mathcal{O}(1/\epsilon)$ & $\mathcal{O}(1/\epsilon)$ & $\cO(1/\epsilon)$ & \textbf{yes} & \textbf{yes}\\ \hline
    \end{tabular}
  }
  \vspace{-0.2cm}
  \caption{Comparison of our method with the previous work}
  \label{tab:prev}
  \vspace{-0.3cm}
\end{table*}

\begin{mthm}
  Let $\{\bx^{(k)}\}_{k\in\integers_+}$
  denote the sequence of \fal~iterates. 
  Then 
$F^*= \lim_{k\in\integers_+}\sum_{i\in\cN}F_i(x^{(k)}_i)$. Furthermore,
$\bx^{(k)}$ is $\epsilon$-optimal and $\epsilon$-feasible
within $\cO(\log(\epsilon^{-1}))$ \fal\ iterations, requiring $\cO\big(\frac{\psi_{\max}^{1.5}}{d_{\min}}\epsilon^{-1}\big)$
  communications per node, and $\cO(\epsilon^{-1})$ gradient and proximal map computations for $\gamma_i$ and
  $\rho_i$, respectively, 
  where $\psi_{\max}$ denotes the largest eigenvalue of
  the Laplacian 
  of $\cG$, and $d_{\min}$ denotes the minimum
  degree over all nodes. 
\end{mthm}
\vspace{-0.3cm}
\subsection{Previous work}
Given the importance of \eqref{eq:dist_opt}, a number
of different
distributed optimization algorithms have been proposed to solve
\eqref{eq:dist_opt}.
\citet{Duchi12} proposed a dual averaging algorithm to
solve~\eqref{eq:centr} in a distributed fashion over $\cG$ 
when each $F_i$ is convex.
This algorithm computes $\epsilon$-optimal solution 
in
$\mathcal{O}(1/\epsilon^2)$ iterations; however, they do not provide
any guarantees on the consensus violation
$\max\{\norm{\bar{x}_i-\bar{x}_j}_2:\
(i,j)\in\cE\}$. 
\citet{nedic2009distributed}
developed a subgradient method with constant step size
$\alpha>0$ for distributed minimization of \eqref{eq:centr} where
the network topology is  time-varying. 
Setting $\alpha=\cO(\epsilon)$ in their method guarantees that consensus
violation and suboptimality is $\cO(\epsilon)$ in
$\cO(1/\epsilon^2)$ iterations; however, since the step size is
constant none of the errors are not guaranteed to decrease further.
Wei and Ozdaglar~\yrcite{wei2012_1,wei20131}, and
recently \citet{makhdoumi2014broadcast} proposed an
alternating direction method of multipliers~(ADMM) algorithm that
computes an $\epsilon$-optimal  and $\epsilon$-feasible
solution in $\mathcal{O}(1/\epsilon)$ proximal map evaluations for $F_i$. 
There are several problems where one can compute the
proximal map for $\rho_i$
efficiently; however, computing the proximal map for $F_i=\rho_i+\gamma_i$ is
hard -see Section \ref{sec:numerical} for an example. One can
overcome this limitation of
ADMM by locally splitting variables, i.e., setting
$F_i(x_i,y_i):=\rho_i(x_i)+\gamma_i(y_i)$, and adding a constraint $x_i=y_i$ in
\eqref{eq:dist_opt}. This approach \emph{doubles} local memory requirement; 
in addition, in order for
ADMM to be efficient, 
proximal maps for \emph{both} $\rho_i$ and $\gamma_i$ 
must be efficiently computable.
When each $F_i$ 
is \emph{smooth} and has \emph{bounded} gradients, 
\citet{jakovetic2011fast} developed a fast distributed gradient methods
with $\mathcal{O}(1/\sqrt{\epsilon})$ convergence rate. Note
that for the quadratic loss, which is one of the most commonly used loss functions, the gradient  
is \emph{not} bounded.
\citet{chen2012fast} proposed an inexact
proximal-gradient method for distributed minimization of
\eqref{eq:centr} that is able to compute $\epsilon$-feasible and
$\epsilon$-optimal solution 
in $\cO(\epsilon^{-1/2})$
iterations which require 
$\cO(\epsilon^{-1})$
communications per node over a time-varying network topology when
$F_i=\rho+\gamma_i$, assuming that  
the non-smooth term $\rho$ is the \emph{same} at all nodes,
and $\grad\gamma_i$ is \emph{bounded} for all $i\in\cN$. 
In contrast, $\fal$ proposed in this paper
is able to \emph{asynchronously} compute
an $\epsilon$-optimal $\epsilon$-feasible solution in
$\cO(\epsilon^{-1})$ communications per node, allowing node specific non-smooth functions $\rho_i$, and without
assuming bounded $\grad\gamma_i$ for any $i\in\cN$. 

\citet{ia10} proposed an
efficient first-order augmented Lagrangian~(\textbf{FAL}) algorithm for
the basis pursuit problem $\min_{x \in
  \mathbb{R}^n} \left\{ \|x\|_1: Ax=b \right\} $
to compute an $\epsilon$-optimal and
$\epsilon$-feasible solution to  within
$\mathcal{O}(\kappa^2(A)/\epsilon)$ matrix-vector multiplications,
where $A \in \mathbb{R}^{m \times n}$ such that $\Rank(A)=m$, and
$\kappa(A) := \sigma_{\max}(A) / \sigma_{\min}(A)$ denotes the
condition number of $A$. 
In this work, we extend their \textbf{FAL}~algorithm to solve a more
general version of
\eqref{eq:dist_opt} in Section~\ref{sec:fal_conv} and
\ref{sec:fal_subproblems},
and establish the \textbf{Main Result} for \eqref{eq:dist_opt} in Section~\ref{sec:laplacian}.
In Section~\ref{sec:fal_async}, we 
propose an \emph{asynchronous} version of \fal. 
It is important to emphasize that \fal~can be easily extended to solve \eqref{eq:dist_opt} when
there are \emph{global} constraints on network resources of the form $Ex-q\in\cK$,
where $\cK$ is a proper cone, and none of the algorithms
discussed above can accommodate such global conic constraints
efficiently. Due to space limitations, we do not discuss this
extension here; however, the analysis would be similar to \cite{Aybat15_1J,Aybat14_1J}.
\vspace{-0.2cm}
\section{Methodology}
\begin{defn}
\label{def:main_def}
\emph{(a)} Let $\Gamma$ be the set of convex functions $\gamma:  \reals^n
\rightarrow \reals$ such that $\grad\gamma$ is Lipschitz continuous
with constant $L_{\gamma}$, and $\gamma(x) \geq \underline{\gamma}$ for all $x\in \reals^n$
for some $\underline{\gamma} \in \reals$. 
\emph{(b)} Let $\cR$ be the set of convex functions
$\rho: \mathbb{R}^n \rightarrow \mathbb{R}$ 
such that subdifferential of $\rho$ is uniformly bounded on $\reals^{n}$,
i.e., there exists $B >0$ such that $\norm{q}_2 \leq B$ for all $q\in \partial
\rho(x)$,  $x\in \mathbb{R}^n$; and 
$\tau \| x \|_2 \leq \rho(x)$ for all $x \in \mathbb{R}^n$ for some $\tau>0$.
\end{defn}
\begin{assump}\label{assump}
For all $i\in\cN$, we assume that $\gamma_i\in\Gamma$ and $\rho_i\in\cR$
with corresponding constants $L_{\gamma_i}, \underline{\gamma_i}, B_i$ and
$\tau_i$.
\end{assump}

Most of the important regularizers and loss functions used in
machine learning and statistics literature lie in $\cR$ and $\Gamma$,
respectively. In particular, 
any norm, 
e.g., $\norm{\cdot}_\alpha$ with
$\alpha\in\{1,2,\infty\}$, group norm (see Section~\ref{sec:numerical}), nuclear
norm, etc., weighted sum of these norms, e.g., sparse group norm (see
Section~\ref{sec:numerical}), all belong to $\cR$. Given $A\in\reals^{m\times n}$ and $b\in\reals^m$,
quadratic-loss $\norm{Ax-b}_2^2$, Huber-loss $\sum_{i=1}^mh(a_i^\T x-b_i)$ (see
Section~\ref{sec:numerical}), logistic-loss
$\sum_{i=1}^m\log\left(1+e^{-b_i a_i^\T x}\right)$, or
fair-loss~\cite{Blatt07_1J} functions all belong to~$\Gamma$.

Throughout the paper, 
we adopt the notation $\bx=(x_i;\bx_{-i})$ with $\bx_{-i}=(x_j)_{j\neq
  i}$ to denote a vector where $x_i$ and $\bx_{-i}$ are treated as
variable and parameter sub-vectors of $\bx$, respectively. Given
$f:\reals^{nN}\rightarrow\reals$,
$\grad_{x_i}f(\bx)\in\reals^n$ denotes the sub-vector of $\grad
f(\bx)\in\reals^{nN}$ corresponding to components of
$x_i\in\reals^n$.\vspace{-0.25cm}
\subsection{APG Algorithm for the Centralized Model}
\label{sec:central}
Consider the centralized version  \eqref{eq:centr} where  all
the functions $F_i$ are available at a
central node, and all computations are carried out at this node.
Suppose
$\{\rho_i\}_{i\in\cN}$ and $\{\gamma_i\}_{i\in\cN}$
satisfy Assumption~\ref{assump}. Let $\rho(x) := \sum_{i=1}^N
\rho_i(x)$ and $\gamma(x) := \sum_{i=1}^N \gamma_i(x)$. Lipschitz
continuity of each $\grad\gamma_i$ with constant $L_{\gamma_i}$
implies that $\nabla \gamma$ is also Lipschitz continuous with
constant $L_{\gamma} = \sum_{i=1}^N L_{\gamma_i}$. When $\textbf{prox}_{
  {\rho/L_{\gamma}}}$ can be computed
efficiently, 
the accelerated proximal gradient~(APG) algorithm proposed
in~\cite{Beck09,Tseng08} 
guarantees that\vspace{-0.2cm}
\begin{equation}\label{apg}
0 \leq F(x^{(\ell)}) -F^* \leq \frac{2 L_{\gamma}}{(\ell+1)^2}
\norm{x^{(0)} - x^*}_2^2,\vspace{-0.2cm}
\end{equation}
where $x^{(0)}$ is the initial iterate and $ x^* \in \argmin_{x \in
  \mathbb{R}^n} F(x)$ --see Corollary 3 in~\cite{Tseng08}, and
Theorem 4.4 in~\cite{Beck09}. Thus, APG can compute an
$\epsilon$-optimal solution to \eqref{eq:centr} within
$\cO(\sqrt{L_{\gamma}}\epsilon^{-\frac{1}{2}})$ iterations.

As discussed above, the centralized APG algorithm cannot be applied when the
nodes are unwilling or unable to communicate the privately known functions
$\{F_i\}_{i \in \cN}$ to a central node. There are many other setting
where one may want to solve
\eqref{eq:centr} as a ``distributed'' problem. For
instance, although $\prox{t\rho_i}$ can be computed efficiently for all
$t>0$ and $i\in\cN$, $\prox{\rho/L_{\gamma}}$ may be hard to compute. 
As an example, consider a problem with $\rho_1(X)=\sum_{i,j}|X_{ij}|$ and
$\rho_2=\sum_{i=1}^{\Rank(X)}\sigma_i(X)$, where $\sigma(X)$ denotes the
vector of singular values for $X \in\reals^{n_1\times
  n_2}$. Here, $\prox{t\rho_i}$ is easy to compute for all $t >0$ and $i\in\{1,2\}$; however,
$\prox{t(\rho_1+\rho_2)}$ is hard to compute. Thus, the ``centralized''
APG algorithm \emph{cannot} be applied. In the rest of this paper, we focus on
\emph{decentralized} algorithms. \vspace{-0.25cm}
\subsection{\fal~Algorithm for the Decentralized Model}
\label{sec:FAL}
Let $\mb{x} = \big(x_1^\top, \ldots, x_N^\top)^\top \in\reals^{nN}$
denotes a vector formed by  concatenating
$\{x_i\}_{i\in\cN}\subset\reals^n$ as a long column vector. Consider the
following optimization problem of the form: \vspace{-0.2cm}
\begin{equation}\label{eq:general_dist}
\bar{F}^* \! :=\!\min_{ \mb{x}\in \mathbb{R}^{nN} } \Big\{
\bar{F}(\bx):=\bar{\rho}(\mb{x}) + \bar{\gamma}(\mb{x}) \text{ s.t. }
 A\bx=b \Big\}, \vspace{-0.2cm}
\end{equation}
where $\bar{\rho}(\bx) := \sum_{i=1}^N \rho_i (x_i)$, $\bar{\gamma}(\bx)
:= \sum_{i=1}^N \gamma_i (x_i)$, and $A \in \mathbb{R}^{m \times nN}$
has $\Rank(A)=m$, i.e., the linear map is \emph{surjective}. In
Section~\ref{sec:laplacian}, we show that the \emph{distributed}
optimization problem in~\eqref{eq:dist_opt} is a special case
of~\eqref{eq:general_dist}, i.e., for all connected $\cG$\footnote{$\cG$
  can contain cycles.}, there exists a \emph{surjective} $A$ such that
\eqref{eq:dist_opt} is equivalent to \eqref{eq:general_dist}. In the rest
of the section, we will use the following notation: Let
$\{A_i\}_{i\in\cN}\subset\reals^{m\times n}$ such that
$A=[A_1,A_2,\ldots,A_N]$; $\bar{L}:=\max_{i\in\cN}L_{\gamma_i}$,
$\bar{\tau}:=\min_{i\in\cN}\tau_i$.

We propose to solve \eqref{eq:general_dist} by \emph{inexactly} solving
the following
sequence of subproblems in a \emph{distributed} manner: \vspace{-0.2cm}
\begin{align}
\xkopt\in\argmin_{ \bx \in \reals^{nN} } \pk(\bx):=\lk \bar{\rho}(\bx)+\fk(\bx), \label{eq:subproblem}\\
\fk(\bx):=\lk\bar{\gamma}( \bx)+\tfrac{1}{2}\norm{A\bx-b-\lk\tk}_2^2, \label{eq:fk}
\end{align}

\vspace{-0.3cm}
\noindent for appropriately chosen sequences of penalty parameters $\{\lk\}$ and
dual variables $\{\tk\}$ such that $\lk\searrow 0$. In
particular, given $\{\ak,\xik\}$ satisfying $\ak\searrow 0$ and $\xik\searrow 0$, the iterate sequence $\{\xk\}$ is constructed such that every
$\xk$ satisfies one of the following conditions: \vspace{-0.2cm}
\begin{equation}
\label{eq:inexact}
\begin{array}{ll}
(a)  &\pk(\xk)-\pk(\xkopt)\leq\ak,\\
(b) &  \exists \gk_i\in\partial_{x_i} \pk(\bx)|_{\bx=\xk}\ \\
  & \hspace*{2pt} \mbox{ s.t. }
  \max_{i\in\cN}\norm{\gk_i}_2\leq\tfrac{\xik}{\sqrt{N}},
\end{array}
\end{equation}

\vspace{-0.5cm}
\noindent $\partial_{x_i} \pk(\bx)|_{\bx=\mb{\bar{x}}}\!:=\!\lk \partial
\rho_i(x_i)|_{x_i=\bar{x}_i}+\grad_{x_i}\fk(\mb{\bar{x}})$. Note that $\grad\fk(\bx)$ is Lipschitz continuous in $\bx\in\reals^{nN}$ with constant $\lk \bar{L}+\sigma^2_{\max}(A)$. Given $\{\bx^{(0)},\lambda^{(0)},\alpha^{(0)},\xi^{(0)}\}$ and $c\in(0,1)$, we
choose the sequence $\{\lk,\ak,\xik,\tk\}$ 
as shown in Fig.~\ref{fig:fal}.

\begin{figure}[thpb]
\centering
 \framebox{\parbox{3.1in}{\footnotesize
 \textbf{Algorithm DFAL} $\left(\lambda^{(1)},\alpha^{(1)},\xi^{(1)}\right)$ \\
 Step $0$: \emph{Set} $ \theta^{(1)}= \mathbf{0}, k=1$\\
 Step $k$: ($k \geq 1$)\\
 \text{ } 1. \emph{Compute} $\bx^{(k)}$ \emph{such that} \eqref{eq:inexact}(a)
   \emph{or} \eqref{eq:inexact}(b) \emph{holds}\\
 \text{ } 2. $\theta^{(k+1)} =  \theta^{(k)} - \frac{A\xk - b }{\lk} $\\
 \text{ } 3. $\lambda^{(k+1)}=c \lk, \quad \alpha^{(k+1)}=c^2~\ak, \quad \xi^{(k+1)}=c^2~\xik$
}}
\caption{\footnotesize First-order Augmented Lagrangian algorithm}
\label{fig:fal}
\end{figure}

In Section~\ref{sec:fal_conv}, we show that \fal\ can compute an
$\epsilon$-optimal and $\epsilon$-feasible $\bx_{\epsilon}$ to ~\eqref{eq:general_dist},
i.e.,  $ \| A {\bx}_{\epsilon} - b  \|_2 \leq \epsilon$ and $ | \bar{F}(
{\bx}_{\epsilon}) - F^*  | \leq \epsilon$, 
in at most
$\cO(\log(1/\epsilon))$ iterations.

Next, in Section~\ref{sec:fal_subproblems}, we show that 
computing an $\epsilon$-optimal, $\epsilon$-feasible solution $x_\epsilon$ requires at most
$\cO\left(\frac{\sigma^3_{\max}(A)}{\min_{i\in\cN}\sigma^2_{\min}(A_i)}\epsilon^{-1}\right)$ floating
point operations. Using this result, in Section~\ref{sec:laplacian} we
establish that \fal~can compute $x_\epsilon$ in a
distributed manner within $\cO(\epsilon^{-1})$ communication steps, 
i.e., the \textbf{Main Result} stated in Section~\ref{sec:intro}. Finally, in
Section~\ref{sec:fal_async} we show how to modify \fal~for an \emph{asynchronous}
computation setting.

\subsubsection{DFAL iteration complexity}
\label{sec:fal_conv}
We first show that $\{\xk\}$ is a bounded sequence, and then argue that
this also implies boundedness of $\{\tk\}$. 
First, we start with a technical lemma that will be used in establishing the
main results of this section.
\begin{lemma}\label{lem:eps-bound}
Let $\bar{\rho}:\reals^{nN}\rightarrow\reals$ be defined as
$\bar{\rho}(\bx)=\sum_{i\in\cN}\rho_i(x_i)$, where $\rho_i \in \cR$ with uniform bound $B_i$ on its subdifferential for all $x\in\reals^n$ and for all $i\in\cN$.
Let $f:\reals^{nN}\rightarrow\reals$ denote a convex function
such that there exist constants $\{L_i\}_{i=1}^N\subset\reals_{++}$
that satisfy 
\vspace{-0.25cm}
\begin{align*}
  f(\by)\leq f(\bar{\by})+\grad
f(\bar{\by})^\T(\by-\bar{\by})+\sum_{i=1}^N\frac{L_i\norm{y_i-\bar{y}_i}_2^2}{2}
\end{align*}

\vspace{-0.45cm}
for all $\by,\bar{\by}\in\reals^{nN}$.
Given $\alpha,\lambda \geq 0$, and $\mb{\bar{x}} \in \reals^{nN}$
such that $\lambda \rho( \mb{\bar{x}} ) + f(\mb{\bar{x}}) - \min_{\bx \in \reals^{nN}}
\{\lambda \rho(\bx) + f(\bx)\} \leq \alpha$, it follows that $\| \grad_{x_i} f(\mb{\bar{x}})\|_2  \leq \sqrt{2L_i \alpha} + \lambda B_i$ for all $i\in\cN$.
\end{lemma}
In Lemma~\ref{lem:lipschitz-fk} we show that function $\fk$ defined in \eqref{eq:fk} 
satisfies the condition given in Lemma~\ref{lem:eps-bound}.
\begin{lemma}\label{lem:lipschitz-fk}
The function $\fk$ in \eqref{eq:subproblem} satisfies the condition in
Lemma~\ref{lem:eps-bound} with the constants $L_i = \Lk_i$, where
$\Lk_i:=\lk L_{\gamma_i}+\sigma_{\max}^2(A)$ for all
$i\in\cN$. 
\end{lemma}
Lemma~\ref{lem:eps-bound} and Lemma~\ref{lem:lipschitz-fk} allow us to bound
$\norm{\theta^{(k+1)}}_2$ in terms of
$\{\norm{\grad_{x_i}\gamma(x_i^{(k)})}_2\}_{i\in\cN}$. We later use this bound in an inductive
argument to establish that  
the sequence
$\{\xk\}$ is bounded. 
\begin{lemma}
\label{lem:theta-bound}
Let $\{\xk\}$ be the \fal~iterate sequence, i.e., at least one of the
conditions in \eqref{eq:inexact} hold for all $k\geq1$. 
Define $\Theta_i^{(k)}:=\max\left\{\sqrt{2 \Lk_i \frac{\ak}{(\lk)^2}},~\frac{1}{\sqrt{N}}\frac{\xik}{\lk}\right\} +
  B_i+\norm{\grad \gamma_i(x^{(k)}_i)}_2$. Then for all $k\geq 1$, we have
\vspace{-0.2cm}
\begin{align*}
\norm{\theta^{(k+1)}}_2 \leq \min_{i\in\cN}\left\{\frac{\Theta_i^{(k)}}{\sigma_{\min}(A_i)}\right\}.
\end{align*}
\end{lemma}
\vspace{-0.25cm}
Theorem~\ref{thm:bounded-x}
establishes that the \fal~iterate sequence $\{\xk\}$ is  bounded whenever
$\{\rho_i,\gamma_i\}_{i\in\cN}$ satisfy Assumption~\ref{assump}; 
therefore, the sequence of dual variables $\{\tk\}$ is bounded according to Lemma~\ref{lem:theta-bound}. 
%
\begin{thm}
\label{thm:bounded-x}
Suppose Assumption~\ref{assump} holds. 
Then there exist constants
$B_x,\ B_\theta,\ \bar{\lambda}>0$ such that
$\max\{\norm{\bx_*^{(k)}}_2,\norm{\xk}_2\}\leq B_x$ and $\norm{\tk}_2\leq
B_\theta$ for all $k\geq 1$, whenever $\lambda^{(1)}$ and $\xi^{(1)}$ are
chosen such that $0<\lambda^{(1)}\leq\bar{\lambda}$ and
$\frac{\xi^{(1)}}{\lambda^{(1)}}<\bar{\tau}$.
\end{thm}
We are now ready to state a key result that will imply the iteration
complexity of \fal.
\begin{thm}
\label{thm:iter-complexity}
Suppose Assumption~\ref{assump} holds and $\lambda^{(1)}$
and $\xi^{(1)}$ are chosen according to Theorem~\ref{thm:bounded-x}. Then the
primal-dual iterate sequence $\{\xk,\tk\}$ generated by \fal~satisfy
\begin{enumerate}[(a)]
\item $\norm{A\xk-b}_2 \leq 2 B_\theta\lk$,
\item $\bar{F}(\xk)-\bar{F}^*\geq -\lk\frac{\left(\norm{\theta^*}_2+B_\theta\right)^2}{2}$
\item $\bar{F}(\xk)-\bar{F}^*\leq\lk\left(\frac{B_\theta^2}{2}+\frac{\max\left\{\alpha^{(1)},~\xi^{(1)}B_x\right\}}
  {\left(\lambda^{(1)}\right)^2}\right)$,
\vspace{-0.25cm}
\end{enumerate}
where $\theta^*$ denotes any optimal dual solution to
\eqref{eq:general_dist}.
\end{thm}
\vspace{-1mm}
\begin{cor}
\label{cor:N_eps}
The \fal~iterates $\xk$ are $\epsilon$-feasible,
  i.e., $\norm{A\xk-b}_2\leq\epsilon$, and $\epsilon$-optimal, i.e.,
  $|\bar{F}(\xk)-\bar{F}^*|\leq\epsilon$, for all $k\geq N(\epsilon)$ and
  $N(\epsilon)=\log_{\frac{1}{c}}(\frac{\bar{C}}{\epsilon})$ for some $\bar{C}>0$.
\end{cor}
\subsubsection{Overall computational complexity for the synchronous algorithm}
\label{sec:fal_subproblems}
Efficiency of \fal~depends on the complexity of the oracle for
Step~1 in Fig.~\ref{fig:fal}. In this section, we construct an oracle \apg~
that computes an $\xk$ satisfying
\eqref{eq:inexact} within $\cO(1/\lk)$ gradient and prox
computations. This result together with Theorem~\ref{thm:iter-complexity}
guarantees that for any $\epsilon>0$, \fal~can compute an
$\epsilon$-optimal and $\epsilon$-feasible iterate within
$\cO\left(\epsilon^{-1}\right)$ floating point operations. Following lemma gives the iteration
complexity of the oracle \apg~displayed in Fig.~\ref{fig:mapg}. 
\begin{lemma}
\label{lem:msapg}
Let $\bar{\rho}:\reals^{nN}\rightarrow\reals$ such that
$\bar{\rho}(\bx)=\sum_{i\in\cN}\rho_i(x_i)$, where
$\rho_i:\reals^n\rightarrow\reals$ is a convex function for all $i\in\cN$,
and $f:\reals^{nN}\rightarrow\reals$ be a convex function such that it satisfies the condition in
Lemma~\ref{lem:eps-bound} for
some constants $\{L_i\}_{i=1}^N\subset\reals_{++}$.
Suppose that
$\by^*\in\argmin\Phi(\by):=\bar{\rho}(\by)+f(\by)$. Then the \apg~iterate sequence
$\{\by^{(\ell)}\}_{\ell\in\integers_+}$, computed as in Fig.~\ref{fig:mapg},
satisfies \vspace{-0.5cm}
\begin{equation}
\label{eq:msapg-rate}
0\leq
\Phi(\by^{(\ell)})-\min_{\by\in\reals^{nN}}\Phi(\by)
\leq \frac{ \sum_{i=1}^N 2 L_i \norm{y_i^{(0)} - y_i^*}_2^2}{(\ell+1)^2}.
\end{equation}
\end{lemma}
\vspace{-0.4cm}
\proof
\eqref{eq:msapg-rate} follows from adapting the proof of Theorem~4.4
in~\citet{Beck09} for the case here.
\endproof
\vspace{-0.25cm}
\begin{figure}[h!]
\centering
 \framebox{\parbox{3.1in}{\footnotesize
 \textbf{Algorithm MS-APG} ( $\bar{\rho}, f, \by^{(0)}$ ) \\[1.5mm]
 Step $0$: Take $ \bar{\by}^{(1)} = \by^{(0)}, t^{(1)}=1 $\\
 Step $\ell$: ($\ell \geq 1$)\\
 \text{ } 1. $\by_i^{(\ell)} =  \textbf{prox}_{\rho_i/L_i} \left(
     \bar{\by}_i^{(\ell)} - \nabla_{y_i} f( \bar{\by}^{(\ell)} )/ { L_i }
   \right)\quad \forall i\in\cN$\\
 \text{ } 2. $t^{(\ell+1)}= ( 1+ \sqrt{ 1+ 4\left(t^{(\ell)}\right)^2} )/2$\\
 \text{ } 3. $\bar{\by}^{(\ell+1)}
 = \by^{(\ell)} + \frac{t^{(\ell)} -1}{t^{(\ell+1)}} \left( \by^{(\ell)} -
   \by^{(\ell-1)} \right )$
}}
\vspace{-0.05cm}
\caption{\footnotesize Multi Step - Accelerated Prox. Gradient~(MS-APG) alg.}
\label{fig:mapg}
\end{figure}

Consider the problem $\Phi^*=\min\Phi(\by):=\bar{\rho}(\by)+f(\by)$ defined in
Lemma~\ref{lem:msapg}. Note that $\grad f$ is Lipschitz continuous with
constant $L=\max_{i\in\cN}L_i$. In \apg~algorithm, the step length $1/L_i\geq 1/L$
is different for each $i\in\cN$.
Instead, 
if one were to use the APG algorithm~\cite{Beck09,Tseng08}, 
then the step length  
would have been $1/L$ for all $i \in \cN$. 
When $\{L_i\}_{i\in\cN}$ are close to each other, the performances of \apg~and
APG are on par; however, when
$\frac{\max_{i\in\cN}L_i}{\min_{i\in\cN}L_i} \gg 1$, APG 
can only take very tiny steps for all $i\in\cN$;  hence, 
\apg\ is likely to converge much faster in practice.

Since the subproblem \eqref{eq:subproblem} is in the form given in
Lemma~\ref{lem:msapg}, the following result  
immediately follows. \vspace{-0.25cm}
\begin{lemma}
\label{lem:subiter-complexity}
The iterate sequence $\{\by^{(\ell)}\}_{\ell\in\integers_+}$ generated
when we call \apg$\left(\lk\bar{\rho},\fk,\bx^{(k-1)}\right)$ satisfies
$\pk(\by^{(\ell)})-\pk(\xkopt)\leq\ak$, for all $\ell\geq
\sqrt{\frac{ \sum_{i=1}^N 2\Lk_i \norm{x_i^{(k-1)} -
      x^{(k)}_{\ast i}}_2^2}{\ak}}-1$,
where $\Lk_i$ is defined in Lemma~\ref{lem:lipschitz-fk} and
$x_{*i}^{(k)}$ represents the $i$-th block of $x_*^{(k)}$. Hence, one can
compute $\xk$ satisfying \eqref{eq:inexact} within $\cO(1/\lk)$
\apg~iterations.
\vspace{-0.3cm}
\end{lemma}
Theorem~\ref{thm:iter-complexity} and Lemma~\ref{lem:subiter-complexity}
together imply that \fal~can compute an $\epsilon$-feasible, and
$\epsilon$-optimal solution to \eqref{eq:general_dist} within
$\cO(1/\epsilon)$ \apg~iterations. Due to space considerations, we will
only state and prove this result for the case where $\grad\bar{\gamma}$ is
bounded in $\reals^{nN}$ since the bounds $B_\theta$ and $B_x$ are more
simple for this case. Note that Huber-loss, logistic-loss,
and fair-loss functions indeed have bounded gradients.
\begin{thm}
\label{thm:rate}
Suppose that $\exists G_i>0$ such that $\norm{\grad \gamma_i(x)}_2\leq
G_i$ for all $x\in\reals^n$ and for all $i\in\cN$. Let
$N^{o}_{\fal}(\epsilon)$ and $N^{f}_{\fal}(\epsilon)$ denote the number of
\fal-iterations to compute an $\epsilon$-optimal, and an
$\epsilon$-feasible solutions to \eqref{eq:general_dist}, respectively.
Let $N^{(k)}$ denote \apg~iteration number required to
compute $\xk$ satisfying at least one of the conditions in \eqref{eq:inexact}. Then \vspace{-0.2cm}
\begin{equation*}
\begin{array}{c}
\sum_{k=1}^{N^{o}_{\fal}(\epsilon)}N^{(k)}
=\cO\left(\Theta^2\sigma_{\max}(A)\epsilon^{-1}\right),\\[0.2cm]
\sum_{k=1}^{N^{f}_{\fal}(\epsilon)}N^{(k)}=\cO\left(\Theta\sigma_{\max}(A)\epsilon^{-1}\right),
\end{array}
\end{equation*}

\vspace{-0.4cm}
\noindent where $\Theta=\frac{\sigma_{\max}(A)}{\min_{i\in\cN}\sigma_{\min}(A_i)}$. \vspace{-0.25cm}
\end{thm}
\subsubsection{Synchronous Algorithm for distributed optimization}
\label{sec:laplacian}
In this section, we show that the decentralized optimization
problem~\eqref{eq:dist_opt} is a special case of~\eqref{eq:general_dist};
therefore, Theorem~\ref{thm:rate} establishes the
\textbf{Main Result} stated in the Introduction. We also show that the steps in \fal\ can
be further simplified in this context.

Construct a directed graph by introducing an arc $(i,j)$ where $i < j$ for
every edge $(i,j)$ in the undirected graph $\cG = (\cN,\cE)$. 
Then the constraints $x_i - x_j=0$ for all $(i, j) \in \mathcal{E}$ in the
distributed optimization problem~\eqref{eq:dist_opt}
can be reformulated as
$C\bx=0$, where $C\in\reals^{n|\cE|\times nN}$ is a block matrix such that
the block $C_{(i, j), l} \in \mathbb{R}^{n \times n}$ corresponding to
the edge $(i,j)\in\cE$ and node $l\in\cN$, i.e.,
$C_{(i, j), l}$ is equal to $I_n$ if $l=i$, $-I_n$ if $l=j$, and $\mb{0}_n$ otherwise,
where $I_n$ and $\mb{0}_n$ denote $n\times n$ identity and zero matrices,
respectively. Let $\Omega \in \reals^{N \times N}$ be the Laplacian of
$\cG$, i.e., for all $i\in\cN$, $\Omega_{ii}=d_i$, and for all
$(i,j)\in\cN\times\cN$ such that $i\neq j$, $\Omega_{ij}=-1$ if either
$(i, j)\in\cE$ or $(j,i)\in\cE$, where $d_i$ denotes the degree of
$i\in\cN$. Then it follows that \vspace{-0.3cm}
\[
\Psi := C^{\mathsf{T}}C = 
\Omega \otimes I_n, \vspace{-0.05cm}
\]
where $\otimes$ denotes the Kronecker product.
Let $\psi_{\max} := \psi_1\geq\psi_2\geq\ldots\geq\psi_N$ be the eigenvalues of
$\Omega$. Since $\cG$ is connected, $\Rank(\Omega) = N-1$,
i.e., $\psi_{N-1}>0$ and $\psi_N=0$. From the structure of $\Psi$ it
follows that that $\{\psi_i\}_{i=1}^N$ are also the eigenvalues of
$\Psi$, each with
\emph{algebraic multiplicity} $n$. Hence, $\Rank(C)=n(N-1)$.

Let $C = U\Sigma V^{\mathsf{T}}$ denote the reduced singular value
decomposition~(SVD) of $C$,
where $U \in \mathbb{R}^{n |\mathcal{E} | \times n(N -1)}$, $\Sigma =
\text{diag}(\sigma)$, $\sigma \in \mathbb{R}_{++}^{n(N -1)}$, and $V \in
\mathbb{R}^{n N\times n(N -1)}$. Note that
$\sigma^2_{\max}(C)=\psi_{\max}$, and $\sigma^2_{\min}(C)=\psi_{N-1}$. Define $A := \Sigma V^{\mathsf{T}}$. $A\in\reals^{n(N -1) \times
  nN}$ has linearly independent rows;  more importantly, $A^{\mathsf{T}}A
= C^{\mathsf{T}}C = \Psi$; hence, $\sigma^2_{\max}(A)=\psi_{\max}$, and
$\sigma^2_{\min}(A)=\psi_{N-1}$. We also have $\{
\bx\in\reals^{nN} : A\bx=0  \} = \{\bx\in\reals^{nN} : C\bx=0  \}$. Hence,
the general problem in \eqref{eq:general_dist} with $A := \Sigma
V^{\mathsf{T}}$ and $b=\mb{0}\in\reals^{n(N-1)}$ is equivalent to
\eqref{eq:dist_opt}. 
Let $A_i\in\reals^{n(N-1)\times nN}$ and
$C_i\in\reals^{n|\cE|\times nN}$ be the submatrices of $A$ and $C$,
respectively, corresponding to $x_i$, i.e., $A=[A_1,A_2,\ldots,A_N]$, and
$C=[C_1,C_2,\ldots,C_N]$. Clearly, it follows from the definition of $C$
that $\sigma_{\max}(C_i)=\sigma_{\min}(C_i)=\sqrt{d_i}$ for all
$i\in\cN$. Using the property of SVD, it can also be shown for $A= \Sigma
V^{\mathsf{T}}$ that $\sigma_{\max}(A_i)=\sigma_{\min}(A_i)=\sqrt{d_i}$
for all $i\in\cN$. Thus, Theorem~\ref{thm:rate} establishes the \textbf{Main Result}. 

We now show that we do not have to compute the SVD of $C$, or $A$, or even
the dual multipliers $\tk$ when \fal\ is used to solve \eqref{eq:dist_opt}.
In \fal\ the matrix $A$  is used 
in Step~1 (i.e. within the oracle \apg) to compute $\grad\fk$, 
and in Step~2 to compute $\theta^{(k+1)}$. 
Since
$\theta^{(1)}=\mb{0}$, Step~2 in \fal~and \eqref{eq:fk} imply that
$\theta^{(k+1)}=-\sum_{t=1}^{k}\frac{Ax^{(t)}}{\lambda^{(t)}}$, and $\grad
\fk(\bx)=\lk\grad\bar{\gamma}(\bx)+A^\T(A\bx-\lk\tk)
=\lk\grad\bar{\gamma}(\bx)+\Psi\left(\bx
  +\lk\sum_{t=1}^{k-1}\frac{1}{\lambda^{(t)}}\bx^{(t)}\right)$.
Moreover, from the definition of $\Psi$, it follows that \vspace{-0.1cm}
\begin{align*}
\lefteqn{\grad_{x_i}\fk(\bx)=}\\
& & \lk\grad\gamma_i(x_i) +d_i\left(x_i+\bar{x}_i^{(k)}\right)-\sum_{j\in\cO_i}\left(x_j+\bar{x}_j^{(k)}\right),
\end{align*}

\vspace{-0.4cm}
\noindent where $\bar{\bx}^{(k)}:= \sum_{t=1}^{k-1}\frac{\lk}{\lambda^{(t)}}\bx^{(t)}$, and $\cO_i$ denotes the set of nodes adjacent to $i\in\cN$. Thus, it follows that Step~1 of \apg~can be computed in a distributed manner
by only communicating with the adjacent nodes without explicitly computing $\tk$
in Step~2 of \fal.

In particular, for the $k$-th \fal~iteration, each node
$i\in\cN$ stores $\bar{x}_i^{(k)}$ and $\{\bar{x}_j^{(k)}\}_{j\in\cO_i}$,
which can be easily computed \emph{locally} if $\{x^{(t)}_j\}_{j\in\cO_i}$ is
transmitted to $i$ at the end of Step~1 of the previous \fal~iterations $1\leq
t\leq k-1$. Hence, during the $\ell$-th iteration of
\apg$\left(\lk\bar{\rho},\fk,\bx^{(k-1)}\right)$ call, each node $i\in\cN$
can compute $\grad_{y_i}\fk(\bar{\by}^{\ell})$ locally if
$\{\bar{y}^{(t)}_j\}_{j\in\cO_i}$ is transmitted to $i$ at the end of
Step~3 in \apg. It is important to note that every node can independently
check \eqref{eq:inexact}(b), i.e., $\exists
\gk_i\in\partial\rho_i(x_i)|_{x_i=x_i^{(k)}}+\grad_{x_i}\fk(\xk)$ for all
$i\in\cN$ such that $\max_{i\in\cN}\norm{\gk_i}_2\leq\frac{\xik}{\sqrt{N}}$. Hence,
nodes can reach a consensus to move to the next \fal~iteration without
communicating their private information. If \eqref{eq:inexact}(b) does not
hold for $\ell^{(k)}_{\max}:=B_x\sqrt{\frac{2\sum_{i\in\cN}\Lk_i}{\ak}}$
\apg~iterations, then Lemma~\ref{lem:subiter-complexity} implies that
\eqref{eq:inexact}(a) must be true. Hence, all the nodes move to next
\fal~iteration after $\ell^{(k)}_{\max}$ many \apg~updates. For
implementable version of \fal, see Figure~\ref{fig:dfal}, where $B_x$ is the bound in Theorem~\ref{thm:bounded-x}, and $\cN_i:=\cO_i\cup\{i\}$.

\begin{figure}[!h]
    \rule[0in]{3.1in}{1pt}\\
    \textbf{Algorithm DFAL}$~(\bx^{(0)},\li,\ai,\xii,B_x,\psi_{\max})$\\
    \rule[0.125in]{3.1in}{0.1mm}
    \vspace{-0.25in}
    {\footnotesize
    \begin{algorithmic}[1]
    \STATE $k\gets1,\quad \bar{x}_i^{(1)}\gets \mathbf{0},\quad \forall i\in\cN$
    \WHILE{$k\geq 1$} \label{algeq:stop_ALCC}
    \STATE $\ell\gets 1,\quad t^{(1)}\gets 1,\quad
    \mathrm{STOP}\gets\mathbf{false}$
    \STATE $y_i^{(0)}\gets x_i^{(k-1)},\quad \bar{y}_i^{(1)}\gets
    x_i^{(k-1)},\quad \forall i\in\cN$
    \STATE $\Lk_i\gets\lk L_{\gamma_i}+\psi_{\max},\ \forall i\in\cN$
    \STATE $\ell^{(k)}_{\max}\gets B_x\sqrt{\frac{2\sum_{i\in\cN}\Lk_i}{\alpha_k}}$
    \WHILE{$\mathrm{STOP}=\mathbf{false}$} \label{algeq:stop}
    \FOR{$i \in N$}
        \STATE $q_i^{(\ell)}\gets\lk\grad\gamma_i\left(\bar{y}_i^{(\ell)}\right)+
        \sum_{j\in \cN_i} \Omega_{ij}
    \left(\bar{y}_j^{(\ell)}+\bar{x}_j^{(k)}\right)$
    \STATE $y_i^{(\ell)} \gets  \textbf{prox}_{\lk\rho_i/\Lk_i} \left(
      \bar{\by}_i^{(\ell)} - q_i^{(\ell)}/ \Lk_i\right) $ 
    \ENDFOR
    \IF {$\exists g_i\in q_i^{(\ell)}+\lk\partial\rho_i\left(\bar{y}_i^{(\ell)}\right)$ s.t. $\max\limits_{i\in\cN}\norm{g_i}_2\leq\tfrac{\xik}{\sqrt{N}}$}
        \STATE $\mathrm{STOP}\gets\mathbf{true},\quad x_i^{(k)}\gets \bar{y}_i^{(\ell)},\ \forall i\in\cN$
    \ELSIF{$\ell=\ell^{(k)}_{\max}$}
        \STATE $\mathrm{STOP}\gets\mathbf{true},\quad x_i^{(k)}\gets y_i^{(\ell)},\ \forall i\in\cN$
    \ENDIF
    \STATE $t^{(\ell+1)}\gets ( 1+ \sqrt{ 1+ 4\left(t^{(\ell)}\right)^2} )/2$
    \STATE $\bar{y}_i^{(\ell+1)}\gets y_i^{(\ell)} + \frac{t^{(\ell)} -1}{t^{(\ell+1)}} \left(y_i^{(\ell)}-y_i^{(\ell-1)} \right ),\ \forall i\in\cN$
    \STATE $\ell \gets \ell + 1$
    \ENDWHILE
    \STATE $\lambda^{(k+1)}\gets c\lk,\ \alpha^{(k+1)}\gets c^2\ak,\ \xi^{(k+1)}\gets c^2\xik$
    \STATE $\bar{x}_i^{(k+1)}\gets \frac{\lambda^{(k+1)}}{\lk}\left(\bar{x}_i^{(k)}+x_i^{(k)}\right),\quad \forall i\in\cN$
    \STATE $k \gets k + 1$
    \ENDWHILE
    \end{algorithmic}
    \rule[0.25in]{3.1in}{0.1mm}
    }
    \vspace{-0.75cm}
    \caption{\small Dist. First-order Aug. Lagrangian~(DFAL) alg.}
    \label{fig:dfal}
\end{figure}

\vspace{-3mm}
\subsubsection{Asynchronous implementation}
\label{sec:fal_async}
Here we propose an \emph{asynchronous} version of \fal.
Due to limited space, and for the sake of simplicity of the
exposition, we only consider a simple randomized block coordinate
descent~(RBCD) method, which will lead to an \emph{asynchronous}
implementation of \fal~that can compute an $\epsilon$-optimal and
$\epsilon$-feasible solution to~\eqref{eq:dist_opt} with probability $1-p$
within $\cO\left(\frac{1}{\epsilon^2}\log\left(\frac{1}{p}\right)\right)$
\bcd~iterations. In Section~\ref{sec:abcd} of the appendix, we
discuss how to improve this rate to
$\cO\left(\frac{1}{\epsilon}\log\left(\frac{1}{p}\right)\right)$ using an
\emph{accelerated} \bcd. \vspace{-0.25cm}

\begin{figure}[h!]
\centering
 \framebox{\parbox{3.1in}{\footnotesize
 \textbf{Algorithm RBCD} ( $\bar{\rho}, f, \by^{(0)}$ ) \\[1.5mm]
 Step $\ell$: ($\ell \geq 0$)\\
 \text{ } 1. $i\in\cN$ is realized with probability $\frac{1}{N}$\\
 \text{ } 2. $\by_i^{(\ell+1)} =  \prox{\rho_i/L_i} \left(y_i^{(\ell)} - \nabla_{y_i} f( \by^{(\ell)} )/ { L_i }\right)$ \\
 \text{ } 3. $\by_{-i}^{(\ell+1)} = \by_{-i}^{(\ell)}$
}}
\caption{\small Randomized Block Coordinate Descent~(RBCD) alg.}
\label{fig:rbcd}
\vspace{-0.1cm}
\end{figure}

Nesterov~\yrcite{Nesterov12} proposed an RBCD method for solving $\min_{\by\in\reals^{nN}}f(\by)$, where $f$
is convex with \emph{block Lipschitz continuous gradient}, i.e., $
\nabla_{y_i} f(y_i;\by_{-i}) $ is Lipschitz continuous in $y_i$ with
constant $L_i$ for all $i$. Later, \citet{Richtarik12_1J} extended the convergence rate results to 
$\min_{\by \in \mathbb{R}^{nN}} \Phi(\by) := \sum_{i=1}^N \rho_i(y_i) +
f(\by)$, 
such that $\prox{t\rho_i}$ can be computed
efficiently for all $t>0$ and $i\in\cN$, and established that given
$\alpha>0$, and $p\in(0,1)$, 
for $\ell\geq \frac{ 2N C }{\alpha} \left(  1+ \log \frac{1}{p}
\right)$, the
iterate sequence $\{\by^{(\ell)}\}$
computed by \bcd~displayed in Fig.~\ref{fig:rbcd} satisfies \vspace{-0.15cm}
\begin{align}
\label{eq:rbcd_complexity}
\mathbb{P}( \Phi(\by^{(\ell)}) - \Phi^* \leq \alpha) \geq 1- p,
\end{align}

\vspace{-0.35cm}
\noindent where $C:=\max\{ \mathcal{R}_L^2(\by^{(0)}),
  \Phi(\by^{(0)}) - \Phi^* \}$, 
$\mathcal{R}_L^2(\by^{(0)}) := \max\limits_{\by,\by^*}\big\{\sum_{i=1}^N
L_i\norm{y_i-y^*_i}_2^2:\ \Phi(\by)\leq\Phi(\by^{(0)}),\ \by^*\in
\cY^*\big\}$, and $\cY^*$ denotes the set of optimal solutions. \bcd~is
significantly faster in practice for \emph{very large scale} problems,
particularly when the partial gradient $\nabla_{y_i} f(\by)$ can be
computed more efficiently as compared to the full gradient $\nabla
f(\by)$. The \bcd\ algorithm can be implemented for the distributed minimization
problem when the nodes in $\cG$ work
\emph{asynchronously}. Assume that for
any $\by = (y_i)_{i \in \cN}\in\reals^{nN}$, each node  $i$
is \emph{equally likely} to be the first to complete computing $\prox{\rho_i/L_i}
\left(y_i - \grad_{y_i} f( \by )/ L_i\right)$, i.e., each node has an exponential clock with equal rates.
Suppose node $i\in\cN$ is the first node to complete Step~2 of
\bcd. 
Then, instead of waiting
for the other nodes to finish, node $i$ sends a message to its neighbors $j\in\cO_i$
to terminate their computations, and
shares $y^{(\ell+1)}_i$ with 
them. Note that \bcd~can be easily incorporated into \fal~as an
oracle to solve subproblems in~\eqref{eq:subproblem} by replacing
\eqref{eq:inexact}(a) with \vspace{-0.15cm}
\begin{equation}
\label{eq:inexact_asynch}
\mathbb{P} \left( P^{(k)} \big( \mathbf{x}^{(k)} \big) -  P^{(k)} \big(
  \mathbf{x}_*^{(k)} \big) \leq \ak \right) \geq \big(1 - p
\big)^{\frac{1}{N(\epsilon)}},
\end{equation}
where $N(\epsilon) = \log_{\frac{1}{c}}\left(\frac{\bar{C}}{\epsilon}\right)$
defined in Corollary~\ref{cor:N_eps}. Since $\big(1 - p
\big)^{\frac{1}{N(\epsilon)}} \leq 1 - \frac{p}{N(\epsilon)}$ for $p\in(0,1)$, the total number of \bcd~iterarions for the $k$-th subproblem is bounded: $\Nk \leq \cO\big( \frac{1}{\ak}
  \log\big(\frac{N(\epsilon)}{p}\big)\big) = \cO\big(
  \frac{1}{\ak} \big(\log\big(\frac{1}{p}\big) + \log
    \log \big(\frac{1}{\epsilon}\big)\big)\big)$. Hence,
Corollary~\ref{cor:N_eps} and \eqref{eq:rbcd_complexity} imply that
\emph{asynchronous} \fal, i.e., \eqref{eq:inexact}(a) replaced with
\eqref{eq:inexact_asynch}, can compute an $\epsilon$-optimal and
$\epsilon$-feasible solution to~\eqref{eq:dist_opt} with probability $1-p$
within $\cO\left(\frac{1}{\epsilon^2}\log\left(\frac{1}{p}\right)\right)$
\bcd~iterations. These results can be extended to the case where each node has \emph{different} clock rates using~\cite{Qu14_1J}.\vspace{-0.25cm}
\section{Numerical results}
\label{sec:numerical}
\vspace{-0.1cm}
  In this section, we compared \fal~with an ADMM method
  proposed in~\cite{makhdoumi2014broadcast} on the \emph{sparse group LASSO} problem with
  Huber loss: \vspace{-.25cm}
  \begin{equation}
    \hspace*{-2pt}
    \min_{x \in \reals^n} \sum_{i=1}^N \beta_1 \norm{x}_1 + \beta_2 \norm{x}_{G_i}+\sum_{j=1}^{m_i}
    h_{\delta}\left( \mb{a}^\T_i(j)x- b_i \right), \label{eq:problem}\vspace{-0.25cm}
  \end{equation}
  where $\beta_1, \beta_2>0$, $A_i\in\reals^{m_i\times n}$,
  $b_i\in\reals^{m_i}$, $\mb{a}^\T_i(j)$ denotes the $j$-th row of
  $A_i$, the Huber loss function $h_{\delta}(x) := \max\{tx-t^2/2: t \in
  [-\delta,\delta]\}$, and $\norm{x}_{G_i} :=
  \sum_{k=1}^{K}\norm{x_{g_i(k)}}_2$
  denotes the group norm with respect to the partition $G_i$ 
  of $[1,n]:=\{1,\cdots,n\}$ for all $i\in\cN$,
  i.e., $G_i=\{g_i(k)\}_{k=1}^{K}$ such that
  $\bigcup_{k=1}^{K}g_i(k)=[1,n]$, and $g_i(j)\cap g_i(k) =\emptyset$
  for all $j \neq k$.
  In this case,
  $\gamma_i(x):=\sum_{j=1}^{m_i} h_{\delta}\left( \mb{a}^\T_i(j)x- b_i
  \right)$ and $\rho_i(x):=\beta_1\norm{x}_1+\beta_2\norm{x}_{G_i}$.
  Next, we briefly describe the ADMM algorithm in~\cite{makhdoumi2014broadcast}, and
  propose a more efficent variant, \sadmm, for \eqref{eq:problem}. 
  \vspace{-0.1cm}
  \begin{figure}[h!]
\centering
 \framebox{\parbox{3.2in}{
 {\footnotesize
 \textbf{Algorithm SADMM} ( $c,\bx^{(0)}$ ) \\[1.5mm]
 Initialization: $\by^{(0)}=\bx^{(0)}$,  $p_i^{(k)}=\tilde{p}_i^{(k)}=\mathbf{0},\ \ i\in\cN$\\
 Step $\ell$: ($\ell \geq 0$) For $i\in\cN$ compute\\
 \text{ } 1. $x_i^{(k+1)}=\prox{\frac{1}{c(d_i^2+d_i+1)}\rho_i}\left(\tilde{x}_i^{(k)}\right)$\\
 \text{ } 2. $y_i^{(k+1)}=\prox{\frac{1}{c(d_i^2+d_i+1)}\gamma_i}\left(\tilde{y}_i^{(k)}\right)$\\
 \text{ } 3. $s_i^{(k+1)}=\sum_{j\in\cN_i}\Omega_{ij}x_j^{(k+1)}/(d_i+1)$\\
 \text{ } 4. $p_i^{(k+1)}=p_i^{(k)}+s_i^{(k+1)}$\\
 \text{ } 5. $\tilde{s}_i^{(k+1)}=\sum_{j\in\cN_i}\Omega_{ij}y_j^{(k+1)}/(d_i+1)$\\
 \text{ } 6. $\tilde{p}_i^{(k+1)}=\tilde{p}_i^{(k)}+\tilde{s}_i^{(k+1)}$
 }}}
\caption{\small Split ADMM algorithm}
\vspace{-0.5cm}
\label{fig:sadmm}
\end{figure}
\begin{table*}[t!]
\centering
\caption{Comparison of \fal, \textbf{AFAL} (Asynchronous \fal), \admm, and \sadmm.\ (Termination time T=1800 sec.)}
{\scriptsize
\begin{tabular}{cl|cc|cc|cc|cc}
\toprule
Size  & Alg. & \multicolumn{2}{|c|}{Rel. Suboptimality}  & \multicolumn{2}{c|}{Consensus Violation~(CV)} & \multicolumn{2}{c|}{CPU Time (sec.)} & \multicolumn{2}{c}{Iterations} \\
\midrule
  & & Case 1 & Case 2 & Case 1 & Case 2 & Case 1 & Case 2 & Case 1 & Case 2 \\ \midrule
&SDPT3~(C)&0&0&0&0&28&85&24&22\\
&APG~(C)&1E-3&N/A&0&N/A&10&N/A&2173&N/A\\
$n_g=100$
&DFAL~(D)&6E-4, 7E-4&4E-4, 6E-4&5E-5, 2E-5&5E-5, 3E-5&\textbf{5, 5}&\textbf{5, 5}&1103, 1022&1105, 1108\\
$N=5$
&AFAL~(D)&3E-4, 8E-4&3E-4, 6E-4&3E-6, 5E-6&2E-6, 5E-6&\textbf{16, 11}&\textbf{17, 11}&9232, 9083&9676, 9844\\
&ADMM~(D)&6E-5, 5E-5&7E-5, 7E-5&1E-4, 1E-4&1E-4, 1E-4&1125, 808&1090, 771&353, 253&363, 261\\
&SADMM~(D)&1E-4, 3E-4&1E-4, 3E-4&1E-4, 1E-4&1E-4, 1E-4&771, 784&772, 804&592, 606&593, 623\\
\midrule
&SDPT3~(C)&0&0&0&0&28&89&24&22\\
&APG~(C)&1E-3&N/A&0&N/A&10&N/A&2173&N/A\\
$n_g=100$
&DFAL~(D)&4E-4, 7E-4&4E-4, 3E-4&6E-5, 1E-5&6E-5, 2E-5&\textbf{14, 12}&\textbf{14, 13}&1794, 1439&1812, 1560\\
$N=10$
&AFAL~(D)&4E-4, 8E-4&2E-4, 8E-4&7E-6, 2E-6&8E-6, 2E-6&\textbf{48, 65}&\textbf{49, 62}&20711, 41125&21519, 41494\\
&ADMM~(D)&9E-3, 2E-2&8E-3, 2E-2&9E-4, 5E-4&8E-4, 4E-4&T, T&T, T&354, 353&373, 372\\
&SADMM~(D)&3E-4, 9E-3&4E-4, 1E-2&2E-4, 1E-3&2E-4, 1E-3&T, T&T, T&867, 883&865, 879\\
\midrule
&SDPT3~(C)&0&0&0&0&806&1653&26&29\\
&APG~(C)&1E-3&N/A&0&N/A&253&N/A&8663&N/A\\
$n_g=300$
&DFAL~(D)&2E-4, 5E-4&2E-4, 4E-4&6E-5, 4E-5&5E-5, 5E-5&\textbf{77, 64}&\textbf{80, 65}&1818, 1511&1897, 1535\\
$N=5$
&AFAL~(D)&1E-4, 6E-4&2E-5, 6E-4&5E-7, 2E-5&6E-8, 2E-5&\textbf{164, 99}&\textbf{273, 99}&21747, 8760&37212, 8736\\
&ADMM~(D)&5E-2, 1E-3&5E-2, 1E-3&5E-3, 1E-3&5E-3, 1E-3&T, T&T, T&109, 118&109, 118\\
&SADMM~(D)&2E-2, 7E-2&2E-2, 8E-2&2E-3, 3E-3&2E-3, 3E-3&T, T&T, T&269, 274&268, 273\\
\midrule
&SDPT3~(C)&0&0&0&0&806&1641&26&29\\
&APG~(C)&1E-3&N/A&0&N/A&253&N/A&8663&N/A\\
$n_g=300$
&DFAL~(D)&1E-4, 6E-4&6E-4, 1E-3&7E-5, 4E-5&9E-5, 5E-5&\textbf{130, 80}&\textbf{122, 82}&2942, 1721&2794, 1769\\
$N=10$
&AFAL~(D)&2E-4, 7E-4&6E-4, 1E-3&8E-7, 1E-5&3E-7, 1E-5&\textbf{350, 294}&\textbf{437, 288}&48214, 29946&63110, 30371\\
&ADMM~(D)&5E-2, 8E-2&5E-2, 8E-2&7E-3, 9E-3&7E-3, 9E-3&T, T&T, T&114, 124&113, 123\\
&SADMM~(D)&3E-1, 3E+0&3E-1, 3E+0&4E-3, 2E-2&4E-3, 2E-2&T, T&T, T&255, 269&256, 268\\
\bottomrule
\end{tabular}
}
\label{tab:1}
\vspace{-5mm}
\end{table*}
  \subsection{A distributed ADMM Algorithm}
  Let
  $\Omega\in\reals^{N\times N}$ denote the Laplacian of the graph $\cG=(\cN,\cE)$,
  $\cO_i$ denote the set of neighboring nodes of
  $i \in \cN$, and define $\cN_i:=\cO_i\cup\{i\}$. Let
  $\cZ_i:=\{z_i\in\reals^{d_i+1}:\
  \sum_{j\in\cN_i}z_{ij}=\mathbf{0}\}$. \citet{makhdoumi2014broadcast} show that
  \eqref{eq:dist_opt} can be equivalently
  written as \vspace{-0.1cm}
  \begin{eqnarray}
    \min_{x_i \in \mathbb{R}^n,z_i\in\cZ_i} &\sum_{i=1}^N
                                              F_i(x_i):=\rho_i(x_i)+\gamma_i(x_i)
                                              \nonumber\\
    \hbox{s.t.}\quad &\Omega_{ij}x_j=z_{ij},\quad i\in\cN,\
                       j\in\cN_i. \label{eq:admm1_cons1}
  \end{eqnarray}

  \vspace{-0.25cm}
  \noindent Only \eqref{eq:admm1_cons1} is penalized when forming the augmented
  Lagrangian, which is alternatingly minimized in $\bx\in\reals^{nN}$,
  $\bz^\T=[z_1^\T,\ldots,z_N^\T]$, where
  $\cZ_i\ni z_i=[z_{ij}]_{j\in\cN_i}\in\reals^{d_i+1}$. 
  \citet{makhdoumi2014broadcast} establish that
  suboptimality and consensus violation converge to $0$ with a rate
  $\cO(1/k)$, and in each iteration every node communicates $3n$ scalars. 
  From now on, we refer to this algorithm 
  that directly works with $F_i$ as \admm. Computing  $\prox{F_i}$ for
  each $i \in \cN$ is the computational bottleneck in each
  iteration of \admm.
  Note that computing $\prox{F_i}$ for \eqref{eq:problem} is almost
  as \emph{hard} as solving the problem. To deal with this issue, we considered
  the following reformulation:
  \begin{eqnarray*}
  \min_{\substack{x_i,y_i \in \mathbb{R}^n,\\
        z_i,\tilde{z}_i\in\cZ_i}}
  &\sum_{i\in\cN}\rho_i(x_i)+\gamma_i(y_i)\\
  \hbox{s.t.}\quad
  &\Omega_{ij}x_j=z_{ij},\quad i\in\cN,\ j\in\cN_i\\
  &\Omega_{ij}y_j=\tilde{z}_{ij},\quad i\in\cN,\ j\in\cN_i\\
  &x_i=q_i,\ y_i=q_i,\quad i\in\cN.
  \vspace{-0.25cm}
  \end{eqnarray*}
  ADMM algorithm for this formulation is displayed in
  Fig.~\ref{fig:sadmm}, where $c>0$ denotes the penalty parameter. Steps
  of \sadmm~can be derived by minimizing the augmented Lagrangian
  alternatingly in $(\bx,\by)$, and in $(\bz,\mathbf{\tilde{z}},\bq)$
  while fixing the other. As in~\cite{makhdoumi2014broadcast},
  computing $(\bz,\mathbf{\tilde{z}},\bq)$ can be avoided by exploiting the
  structure of optimality conditions. Prox centers in \sadmm~are
  \begin{equation*}
  \begin{array}{l}
  \tilde{x}_i^{(k)}=x_i^{(k)}-\frac{\sum_{j\in\cN_i}\Omega_{ji}(s_j^{(k)}+p_j^{(k)})+r_i^{(k)}+(x_i^{(k)}-y_i^{(k)})/2}{d^2_i+d_i+1},\\
  \tilde{y}_i^{(k)}=y_i^{(k)}-\frac{\sum_{j\in\cN_i}\Omega_{ji}(\tilde{s}_j^{(k)}+\tilde{p}_j^{(k)})-r_i^{(k)}-(x_i^{(k)}-y_i^{(k)})/2}{d^2_i+d_i+1},
  \end{array}
  \end{equation*}
  respectively; and $r_i^{(k+1)}=r_i^{(k)}+(x_i^{(k+1)}-y_i^{(k+1)})/2$.
\subsection{Implementation details and numerical results}
The following lemma shows that in \fal~implementation, each node $i\in\cN$
can check \eqref{eq:inexact}(b) very efficiently. For $x \in \reals$, define $\sgn(x)$
as -1, 0 and 1 when $x<0$, $x=0$, and $x>0$, respectively;
and for $x\in\reals^n$, define $\sgn(x)=[\sgn(x_1),\sgn(x_2),\ldots,\sgn(x_n)]^\T$.
\begin{lemma}
\label{lem:subgrad}
Let $f:\reals^n\rightarrow\reals$ be a differentiable function, $G=\{g(k)\}_{k=1}^{K}$ be a partition of $[1,n]$. For $\beta_1,\beta_2>0$, define
$P=\lambda\rho+f$, where $\rho(x):=\beta_1\norm{x}_1+\beta_2\norm{x}_G$.
Then, for all $\bar{x}\in\reals^n$ and $\xi>0$,
there exists $\nu\in\partial
P(x)|_{x=\bar{x}}$ such that $\norm{\nu}_2\leq\xi$ if and only if
$\norm{\pi^*+\omega^*+\grad f(\bar{x})}_2\leq\xi$ for $\pi^*,\omega^*$
such that for each $k$, if $\bar{x}_{g(k)}\neq\mathbf{0}$, then
$\pi^*_{g(k)}=\lambda\beta_1\sgn(\bar{x}_{g(k)})+\left(\mathbf{1}-
  \sgn(|\bar{x}_{g(k)}|)\right)\odot\eta_{g(k)}$,
and
$\omega^*_{g(k)}=\lambda\beta_2~\frac{\bar{x}_{g(k)}}{\norm{\bar{x}_{g(k)}}_2}$;
otherwise, if $\bar{x}_{g(k)}=\mathbf{0}$, then
$\pi^*_{g(k)}=\eta_{g(k)}$, and $\omega^*_{g(k)}$ equals \vspace{-0.2cm}
\begin{equation*}
-\left(\pi^*_{g(k)}+\grad_{x_{g(k)}}f(\bar{x})\right)
\min\left\{1,~\frac{\lambda\beta_2}{\norm{\pi^*_{g(k)}+\grad_{x_{g(k)}} f(\bar{x})}_2}\right\},
\end{equation*}
where $\odot$ denotes componentwise multiplication, and $\eta_{g(k)}=-\sgn\left(\grad_{x_{g(k)}} f(\bar{x})\right)\odot\min\left\{|\grad_{x_{g(k)}} f(\bar{x})|,~\lambda\beta_1\right\}$.
\end{lemma}
Both \fal~and \sadmm~call for $\prox{\rho_i}$. In Lemma~\ref{lem:prox}, we show that it can be computed in closed form. On the other hand, when \admm, and \sadmm~are implemented on \eqref{eq:problem}, one needs to compute $\prox{F_i}$ and $\prox{\gamma_i}$, respectively; however, these proximal operations do not assume closed form solutions. Therefore, in order to be fair, we computed them using an efficient interior point solver MOSEK~(ver. 7.1.0.12).
\begin{lemma}
\label{lem:prox}
Let $\rho(x) = \beta_1 \| x \|_1 + \beta_2 \| x \|_G$. 
For $t>0$ and $\bar{x}\in\reals^{n}$,
$x^p=\prox{t\rho}(\bar{x})$ 
is given by
$
x^{p}_{g(k)} = \eta'_{g(k)} \max \left\{ 1 -
  \frac{t\beta_2}{\norm{\eta_{g(k)}}_2},~ 0  \right\}$, for $1\leq
k\leq K$, where $\eta' = \sgn(\bar{x}) \odot \max \{ |\bar{x}| - t\beta_1, 0 \}$.
\end{lemma}
\vspace{-0.3cm}
In our experiments, the network was either a \emph{star tree}
or a \emph{clique} 
with either $5$ or $10$ nodes.
The remaining problem parameters defining
$\{\rho_i,\gamma_i\}_{i\in\cN}$ were set as follows. We set
$\beta_1=\beta_2=\frac{1}{N}$, $\delta=1$, and $K=10$. Let $n=K n_g$ for
$n_g\in\{100,300\}$, i.e., $n\in\{1000, 3000\}$. We generated partitions
$\{G_i\}_{i\in\cN}$ in two different ways. For test problems in
\textbf{CASE 1}, we created a
single partition $G=\{g(k)\}_{k=1}^K$ by generating $K$ groups
uniformly at random such that $|g(k)| = n_g$ for
all $k$; and set $G_i=G$ for all $i\in\cN$,
i.e., $\rho_i(x)=\rho(x):=\beta_1\norm{x}_1+\beta_2\norm{x}_G$ for all
$i\in\cN$. For the test problems in \textbf{CASE 2}, we created a different partition $G_i$ for each node $i$, in the same
manner as in \textbf{Case 1}.
  For all $i\in\cN$, $m_i = \frac{n}{2N}$, and the elements of
  $A_i\in\reals^{m_i\times n}$ 
  are i.i.d. with standard Gaussian, 
  and we set $b_i = A_i \bar{x}$ for $\bar{x}_j = (-1)^j e^{-(j-1)/n_g}$ for $j\in[1,n]$.

We solved the distributed optimization problem \eqref{eq:dist_opt} using
  \fal, \afal~(\emph{asynchronous} version of \fal~with \emph{accelerated} \bcd~-see the appendix for details), \admm, and \sadmm~for both cases, on both star trees, and cliques, and for
  $N\in\{5,10\}$ and $n_g\in\{100,300\}$. For each problem setting, we randomly generated 5 instances. For benchmarking, we solved the centralized problem~\eqref{eq:centr} using SDPT3 for both cases. Note
  that for \textbf{Case 1}, $\sum_{i\in\cN}\rho_i(x)=\norm{x}_1+\norm{x}_G$
  and its prox mapping can be computed efficiently, while for \textbf{Case 2}, $\sum_{i\in\cN}\rho_i(x)$ does not assume a simple prox map. Therefore, for the first case we were also able to use APG,
  described in Section~\ref{sec:central}, to solve \eqref{eq:centr} by exploiting the result of Lemma~\ref{lem:prox}. All the algorithms are terminated when the relative suboptimality, $|F^{(k)}-F^*|/|F^*|$, is less then $10^{-3}$, and consensus violation, $\mathrm{CV}^{(k)}$, is less than $10^{-4}$, where $F^{(k)}$ equals to $\sum_{i\in\cN}F_i(x_i^{(k)})$ for \fal~and \admm, and to $\sum\limits_{i\in\cN}F_i\left(\frac{x_i^{(k)}+y_i^{(k)}}{2}\right)$ for \sadmm; $\mathrm{CV}^{(k)}$ equals to $\max_{(ij)\in\cE}\norm{x_i^{(k)}-x_j^{(k)}}_2/\sqrt{n}$ for \fal, and \admm, and to $\max\{\max_{(ij)\in\cE}\norm{x_i^{(k)}-x_j^{(k)}}_2,~\max_{i\in\cN}\norm{x_i^{(k)}-y_i^{(k)}}_2\}/\sqrt{n}$ for \sadmm. If the stopping criteria are not satisfied in 30min., we terminated the algorithm and report the statistics corresponding to the iterate at the termination.

In Table~\ref{tab:1}, '\rm{xxx}~(C)' stands for ``algorithm \rm{xxx} is used to solve the \emph{centralized} problem". Similarly, '\rm{xxx}~(D)' for the \emph{decentralized} one. For the results separated by comma, the left and right ones are for the star tree and clique, resp. Table~\ref{tab:1} displays the means over 5 replications for each case.
The number of iterations in each case clearly illustrates the topology of the network plays an important role in the convergence speed of \fal, which coincides to our analysis in Section~\ref{sec:fal_subproblems}.
\newpage
\bibliography{paper}

\begin{thebibliography}{22}
\providecommand{\natexlab}[1]{#1}
\providecommand{\url}[1]{\texttt{#1}}
\expandafter\ifx\csname urlstyle\endcsname\relax
  \providecommand{\doi}[1]{doi: #1}\else
  \providecommand{\doi}{doi: \begingroup \urlstyle{rm}\Url}\fi

\bibitem[Aybat \& Iyengar(2012)Aybat and Iyengar]{ia10}
Aybat, N.~S. and Iyengar, G.
\newblock A first-order augmented lagrangian method for compressed sensing.
\newblock \emph{SIAM Journal on Optimization}, 22\penalty0 (2):\penalty0
  429--459, 2012.

\bibitem[Aybat \& Iyengar(2013)Aybat and Iyengar]{Aybat15_1J}
Aybat, Necdet~Serhat and Iyengar, Garud.
\newblock An augmented lagrangian method for conic convex programming.
\newblock \emph{arXiv preprint arXiv:1302.6322}, 2013.

\bibitem[Aybat \& Iyengar(2014)Aybat and Iyengar]{Aybat14_1J}
Aybat, Necdet~Serhat and Iyengar, Garud.
\newblock A unified approach for minimizing composite norms.
\newblock \emph{Mathematical Programming}, 144\penalty0 (1-2):\penalty0
  181--226, 2014.

\bibitem[Beck \& Teboulle(2009)Beck and Teboulle]{Beck09}
Beck, A. and Teboulle, M.
\newblock A fast iterative shrinkage-thresholding algorithm for linear inverse
  problems.
\newblock \emph{SIAM J. Img. Sci.}, 2\penalty0 (1):\penalty0 183--202, March
  2009.
\newblock ISSN 1936-4954.

\bibitem[Blatt et~al.(2007)Blatt, Hero, and Gauchman]{Blatt07_1J}
Blatt, D., Hero, A., and Gauchman, H.
\newblock A convergent incremental gradient method with a constant step size.
\newblock \emph{SIAM Journal on Optimization}, 18\penalty0 (1):\penalty0
  29--51, 2007.

\bibitem[Chen \& Ozdaglar(2012)Chen and Ozdaglar]{chen2012fast}
Chen, Annie~I and Ozdaglar, Asuman.
\newblock A fast distributed proximal-gradient method.
\newblock In \emph{Communication, Control, and Computing (Allerton), 2012 50th
  Annual Allerton Conference on}, pp.\  601--608. IEEE, 2012.

\bibitem[Duchi et~al.(2012)Duchi, Agarwal, and Wainwright]{Duchi12}
Duchi, J.~C., Agarwal, A., and Wainwright, M.~J.
\newblock Dual averaging for distributed optimization: Convergence analysis and
  network scaling.
\newblock \emph{IEEE Trans. Automat. Contr.}, 57\penalty0 (3):\penalty0
  592--606, 2012.

\bibitem[Fercoq \& Richt{\'a}rik(2013)Fercoq and
  Richt{\'a}rik]{fercoq2013accelerated}
Fercoq, Olivier and Richt{\'a}rik, Peter.
\newblock Accelerated, parallel and proximal coordinate descent.
\newblock \emph{arXiv preprint arXiv:1312.5799}, 2013.

\bibitem[Grone \& Merris(1994)Grone and Merris]{Grone94_1J}
Grone, R. and Merris, R.
\newblock The laplacian spectrum of a graph. ii.
\newblock \emph{SIAM J. Discrete Math.}, 7\penalty0 (2):\penalty0 221--229,
  1994.

\bibitem[Jakovetic et~al.(2011)Jakovetic, Xavier, and Moura]{jakovetic2011fast}
Jakovetic, Dusan, Xavier, Joao, and Moura, J.
\newblock Fast distributed gradient methods.
\newblock 2011.

\bibitem[Lesser et~al.(2003)Lesser, Ortiz~Jr, and Tambe]{lesser2003distributed}
Lesser, Victor, Ortiz~Jr, Charles~L, and Tambe, Milind.
\newblock \emph{Distributed sensor networks: A multiagent perspective},
  volume~9.
\newblock Springer, 2003.

\bibitem[Makhdoumi \& Ozdaglar(2014)Makhdoumi and
  Ozdaglar]{makhdoumi2014broadcast}
Makhdoumi, Ali and Ozdaglar, Asuman.
\newblock Broadcast-based distributed alternating direction method of
  multipliers.
\newblock In \emph{Communication, Control, and Computing (Allerton), 2014 52nd
  Annual Allerton Conference on}, pp.\  270--277. IEEE, 2014.

\bibitem[Mateos et~al.(2010)Mateos, Bazerque, and
  Giannakis]{mateos2010distributed}
Mateos, Gonzalo, Bazerque, Juan~Andr{\'e}s, and Giannakis, Georgios~B.
\newblock Distributed sparse linear regression.
\newblock \emph{Signal Processing, IEEE Transactions on}, 58\penalty0
  (10):\penalty0 5262--5276, 2010.

\bibitem[McDonald et~al.(2010)McDonald, Hall, and
  Mann]{mcdonald2010distributed}
McDonald, Ryan, Hall, Keith, and Mann, Gideon.
\newblock Distributed training strategies for the structured perceptron.
\newblock In \emph{Human Language Technologies: The 2010 Annual Conference of
  the North American Chapter of the Association for Computational Linguistics},
  pp.\  456--464. Association for Computational Linguistics, 2010.

\bibitem[Necoara \& Suykens(2008)Necoara and Suykens]{necoara2008application}
Necoara, Ion and Suykens, Johan~AK.
\newblock Application of a smoothing technique to decomposition in convex
  optimization.
\newblock \emph{Automatic Control, IEEE Transactions on}, 53\penalty0
  (11):\penalty0 2674--2679, 2008.

\bibitem[Nedic \& Ozdaglar(2009)Nedic and Ozdaglar]{nedic2009distributed}
Nedic, Angelia and Ozdaglar, Asuman.
\newblock Distributed subgradient methods for multi-agent optimization.
\newblock \emph{Automatic Control, IEEE Transactions on}, 54\penalty0
  (1):\penalty0 48--61, 2009.

\bibitem[Nesterov(2012)]{Nesterov12}
Nesterov, Y.
\newblock Efficiency of coordinate descent methods on huge-scale optimization
  problems.
\newblock \emph{SIAM Journal on Optimization}, 22\penalty0 (2):\penalty0
  341--362, 2012.

\bibitem[Qu \& Richt{\'a}rik(2014)Qu and Richt{\'a}rik]{Qu14_1J}
Qu, Zheng and Richt{\'a}rik, Peter.
\newblock Coordinate descent with arbitrary sampling i: Algorithms and
  complexity.
\newblock \emph{arXiv preprint arXiv:1412.8060}, 2014.

\bibitem[Richt\'{a}rik \& Tak\'{a}\u{c}(2012)Richt\'{a}rik and
  Tak\'{a}\u{c}]{Richtarik12_1J}
Richt\'{a}rik, P. and Tak\'{a}\u{c}, M.
\newblock Iteration complexity of randomized block-coordinate descent methods
  for minimizing a composite function.
\newblock \emph{forthcoming, Mathematical Programming, Series A}, 2012.
\newblock DOI: 10.1007/s10107-012-0614-z.

\bibitem[Tseng(2008)]{Tseng08}
Tseng, Paul.
\newblock On accelerated proximal gradient methods for convex-concave
  optimization.
\newblock \emph{submitted to SIAM Journal on Optimization}, 2008.

\bibitem[Wei \& Ozdaglar(2012)Wei and Ozdaglar]{wei2012_1}
Wei, Ermin and Ozdaglar, Asuman.
\newblock Distributed alternating direction method of multipliers.
\newblock 2012.

\bibitem[Wei \& Ozdaglar(2013)Wei and Ozdaglar]{wei20131}
Wei, Ermin and Ozdaglar, Asuman.
\newblock On the o (1/k) convergence of asynchronous distributed alternating
  direction method of multipliers.
\newblock \emph{arXiv preprint arXiv:1307.8254}, 2013.

\end{thebibliography}
\bibliographystyle{icml2015}
\onecolumn
\section{Appendix}
\subsection{Proof of Lemma~\ref{lem:eps-bound}}
\proof
Let $\bx \in \mathbb{R}^{nN}$ and $g_i \in \partial \rho_i(x_i)$ for all $i\in\cN$. From convexity of
$\rho_i$ and Cauchy-Schwarz, it follows that
$\rho_i(x_i)\leq\rho(\bar{x}_i)+\norm{g_i}_2 \norm{x_i-\bar{x}_i}_2$ for all $i\in\cN$. Hence, we have
\begin{align*}
\lambda \bar{\rho}(\bx) + f(\bx)
&\leq  \lambda \rho(\mb{\bar{x}})+f(\mb{\bar{x}}) + \sum_{i\in\cN}\left(\lambda B_i \norm{x_i - \bar{x}_i}_2  + \grad_{x_i} f(\mb{\bar{x}})^\T ( x_i - \bar{x}_i )+ \frac{L_i}{2}
\norm{x_i - \bar{x}_i}_2^2\right).
\end{align*}
Minimizing on both sides and using the separability of the right side, we have $\min_{\bx \in \reals^{nN} } \lambda
\bar{\rho}(\bx) + f(\bx) \leq \lambda \bar{\rho}(\mb{\bar{x}}) + f(\mb{\bar{x}}) + \sum_{i\in\cN} \min_{x_i \in
  \reals^{n}}h_i(x_i)$, where $h_i(x_i) := \nabla_{x_i} f(\mb{\bar{x}})^\T ( x_i -
\bar{x}_i ) + \lambda B_i \norm{x_i - \bar{x}_i}_2 + \frac{L_i}{2} \norm{x_i - \bar{x}_i}_2^2 $. Let $\bar{x}^*_i := \argmin_{x_i \in \mathbb{R}^n} h_i(x_i)$.
Then the first-order optimality conditions imply that
$0 \in \grad_{x_i} f(\mb{\bar{x}}) + L_i( \bar{x}^*_i - \bar{x}_i ) + \lambda B_i
\: \partial \norm{x_i - \bar{x}_i}_2 \Big|_{x_i = \bar{x}^*_i}$ for all $i\in\cN$.

Let $\cI:=\{i\in\cN:\ \norm{\grad_{x_i} f(\mb{\bar{x}})}_2 \leq \lambda B_i\}$. For each $i\in\cN$, there are two possibilities.

\textbf{Case 1:} Suppose that $i\in\cI$, i.e., $\norm{\grad_{x_i} f(\mb{\bar{x}})}_2 \leq \lambda
B_i$. 
Since $\min_{x_i\in\reals^n}h_i(x_i)$ has a
unique solution, and $-\grad_{x_i} f(\mb{\bar{x}}) \in \lambda B_i \: \partial \| x_i
- \bar{x}_i \|_2 \Big|_{x_i = \bar{x}_i}$ when $\norm{\grad_{x_i} f(\mb{\bar{x}})}_2 \leq \lambda
B_i$, it follows that $\bar{x}^*_i = \bar{x}_i$ if and only if $\norm{\grad_{x_i} f(\mb{\bar{x}})}_2 \leq \lambda
B_i$. Hence, $h_i(\bar{x}_i^*)=0$.

\textbf{Case 2:} Suppose that $i\in\cI^c:=\cN\setminus\cI$, i.e., $\norm{\grad_{x_i} f(\mb{\bar{x}})}_2 > \lambda
B_i$. 
In this case, $\bar{x}^*_i \neq
\bar{x}_i$. From the first-order optimality condition, we have
$\grad_{x_i} f(\mb{\bar{x}}) + L_i( \bar{x}^*_i - \bar{x}_i ) + \lambda B_i  \frac{\bar{x}^*_i
  - \bar{x}_i}{\| \bar{x}^*_i - \bar{x}_i \|_2 } = 0$.
Let $s_i:=\frac{\bar{x}^*_i - \bar{x}_i}{\| \bar{x}^*_i - \bar{x}_i \|_2 }$ and $t_i
:= \| \bar{x}^*_i - \bar{x}_i \|_2$, then $s_i = \frac{- \grad_{x_i} f(\mb{\bar{x}})}{L_i t_i+\lambda B_i}$. Since $\| s_i \|_2 =1$, it follows that
$t_i = \frac{ \| \grad_{x_i} f(\mb{\bar{x}}) \|_2 - \lambda B_i}{L_i} >0$,  and $s_i=\frac{
  - \grad_{x_i} f(\mb{\bar{x}})}{ \| \grad_{x_i} f(\mb{\bar{x}}) \|_2}$. Hence, $\bar{x}^*_i =
\bar{x}_i - \frac{ \| \grad_{x_i} f(\mb{\bar{x}}) \|_2 - \lambda B_i}{L_i} \frac{ \grad_{x_i} f(\mb{\bar{x}})}{ \| \grad_{x_i} f(\mb{\bar{x}}) \|_2}$, and $h_i(\bar{x}^*_i)=-\frac{(\|
  \grad_{x_i} f(\mb{\bar{x}}) \|_2 - \lambda B_i)^2}{2L_i}$.

From the $\alpha$-optimality of $\mb{\bar{x}}$, it follows that
\begin{align*}
 \sum_{i\in\cI}\frac{(\| \grad_{x_i} f(\mb{\bar{x}}) \|_2 - \lambda B_i)^2}{2L_i} = -\sum_{i\in\cI} h_i(\bar{x}^*_i) \leq \lambda \bar{\rho}(\mb{\bar{x}}) + f(\mb{\bar{x}}) - \min_{\bx \in \mathbb{R}^{nN} } \lambda \bar{\rho}(\bx) + f(\bx) \leq \alpha,
\end{align*}
which implies that $\| \grad_{x_i} f(\mb{\bar{x}})\|_2  \leq \sqrt{2L_i \alpha} + \lambda B_i$ for all $i\in\cI$. Moreover, $\| \grad_{x_i} f(\mb{\bar{x}})\|_2\leq \lambda B_i$ for all $i\in\cI^c$. Hence, the result follows from these two inequalities.
\endproof
\subsection{Proof of Lemma~\ref{lem:lipschitz-fk}}
\proof
For all $i\in\cN$, since $\grad \gamma_i$ is Lipschitz continuous with
constant $L_{\gamma_i}$, for any $\bx,\bar{\bx}\in\reals^{nN}$, we have
$\gamma_i(x_i) \leq \gamma_i(\bar{x}_i) + \grad\gamma_i(\bar{x}_i)^\T (
x_i -\bar{x}_i ) + \frac{L_{\gamma_i}}{2} \norm{x_i -
  \bar{x}_i}_2^2$. Then, it follows that
\begin{align}
\bar{\gamma}( \bx ) &\leq \sum_{i=1}^N \gamma_i(\bar{x}_i) + \grad
\gamma_i(\bar{x}_i)^\T ( x_i -\bar{x}_i ) + \frac{L_{\gamma_i}}{2}
\norm{x_i - \bar{x}_i}_2^2 \nonumber \\
& \leq \bar{\gamma}(\bar{\bx}) + \grad\bar{\gamma}(\bar{\bx})^\T ( \bx -
\bar{\bx} ) + \sum_{i=1}^N \frac{L_{\gamma_i}}{2} \norm{x_i -
  \bar{x}_i}_2^2. \label{eq:lipschitz-gamma}
\end{align}
Let $h^{(k)}(\bx)=\tfrac{1}{2}\norm{A\bx-b-\lk\tk}_2^2$. It follows that
$\grad h^{(k)}$ is Lipschitz continuous with constant
$\sigma_{\max}^2(A)$. Since $\fk=\lk\bar{\gamma}+h^{(k)}$, the result
follows from \eqref{eq:lipschitz-gamma}.
\endproof
\subsection{Proof of Lemma~\ref{lem:theta-bound}}
\proof
Fix $k\geq 1$. Suppose that $\xk$ satisfies \eqref{eq:inexact}(a). Then
Lemma~\ref{lem:eps-bound} implies that for all $i\in\cN$
\[
\norm{\grad_{x_i} \fk(\xk)}_2=\norm{\lk\grad \gamma_i(x^{(k)}_i) +
  A_i^{\mathsf{T}} (A\xk-b-\lk \theta^{(k)} )}_2   \leq \sqrt{2\Lk_i \ak} +
\lk B_i.
\]
Now, suppose that $\xk$ satisfies \eqref{eq:inexact}(b). Then triangular
inequality immediately implies that $\norm{\grad_{x_i}
  \fk(\xk)}_2\leq\xik/\sqrt{N}+\lk B_i$ for all $i\in\cN$. Combining the two inequalities, and
further using triangular Cauchy-Schwarz inequalities, it follows for all $i\in\cN$ that
$\norm{A\xk - b - \lk \tk}_2 \leq \frac{
  \max\left\{\sqrt{2 \Lk_i \alpha_k},~\xik/\sqrt{N}\right\} + \lk \big(B_i + \norm{\grad \gamma(x^{(k)}_i)}_2 \big)}{{\sigma_{\min}(A_i)}}$.
Hence, we conclude by diving the above inequality by $\lk$ and using the definition of $\theta^{(k+1)}$.
\endproof
\subsection{Proof of Theorem~\ref{thm:bounded-x}}
\proof
Let $A=[A_1,A_2,\ldots,A_N]\in\reals^{m\times nN}$ such that $A_i\in\reals^{m\times n}$ for all $i\in\cN$. Throughout the proof we assume that $\sigma_{\max}(A)\geq \sqrt{\max_{i\in\cN}d_i+1}$, and $\sigma_{\min}(A_i)=\sqrt{d_i}\geq 1$ for all $i\in\cN$, where $d_i\geq 1$ is the degree of $i\in\cN$. Indeed, when $A$ is chosen as described in
Section~\ref{sec:laplacian} corresponding to graph $\cG$, recall that we showed
$\sigma^2_{\max}(A)=\psi_1$, where $\psi_1$ is the largest
eigenvalue of the Laplacian $\Omega$ corresponding to $\cG$. It is
shown in~\cite{Grone94_1J} that when $\cG$ is connected, one has $\psi_1\geq\max_{i\in\cN}d_i+1>1$. Hence, $\sigma_{\max}(A)\geq\sqrt{\max_{i\in\cN}d_i+1}>1$. Moreover, for $A$ chosen as described in
Section~\ref{sec:laplacian} corresponding to graph $\cG$, again recall that $\sigma_{\min}(A_i)=\sqrt{d_i}$ for all $i\in\cN$.

To keep notation simple, without loss of generality, we assume that
$\underline{\gamma_i}=0$ for all $i\in\cN$. Hence,
$\bar{\gamma}(\bx)\geq 0$ for all
$\bx\in\reals^{nN}$. 
Let $\bx^*$ be
a minimizer of \eqref{eq:general_dist}. By Lipschitz continuity of $\grad\gamma_i$, we have for all $i\in\cN$
\begin{align}\label{gamma1}
\norm{\grad \gamma(x_i)}_2 
\leq L_{\gamma_i}\norm{x_i - x^*_i}_2 + \norm{\grad\gamma_i(x^*_i)}_2.
\end{align}
We prove the theorem using induction. We show that, for an appropriately
chosen bound $R$, 
$\norm{\xk-\bx^*}_2
\leq R$ implies that $\norm{\bx^{(k+1)} - \bx^\ast}_2 \leq R$, for all $k
\geq 1$.
Fix $k\geq 1$. First, suppose that $\bx^{(k+1)}$ satisfies \eqref{eq:inexact}(a), i.e.
$P^{(k+1)}(\bx^{(k+1)})\leq
P^{(k+1)}(\bx^{*})+\alpha^{(k+1)}$.
By dividing both sides by $\lambda^{(k+1)}$, it follows from
Assumption~\ref{assump}, $A\bx^*=b$, and $f^{(k+1)}(\cdot)\geq 0$ that
\begin{align}
\label{eq:x-bound-alpha}
\bar{\tau}\norm{\bx^{(k+1)}}_2\leq \bar{\rho}(\bx^*)+\bar{\gamma}(\bx^*) +
\frac{\lambda^{(k+1)}}{2} \left(\norm{\theta^{(k+1)}}_2^2 +
  \frac{\alpha^{(k+1)}}{\left(\lambda^{(k+1)}\right)^2}\right).
\end{align}
Next, suppose  $\bx^{(k+1)}$ satisfies \eqref{eq:inexact}(b). 
It follows from convexity of $P^{(k+1)}$ and Cauchy-Schwarz inequality that
$P^{(k+1)}(\bx^{(k+1)})\leq
P^{(k+1)}(\bx^*)+\xi^{(k+1)}\norm{\bx^{(k+1)}-\bx^*}_2$. Again, dividing
both sides by $\lambda^{(k+1)}$, we get
\begin{align}
\label{eq:x-bound-xi}
\bar{\tau}\norm{\bx^{(k+1)}}_2\leq \bar{\rho}(\bx^*)+\bar{\gamma}(\bx^*) +
\frac{\lambda^{(k+1)}}{2} \norm{\theta^{(k+1)}}_2^2 +
\frac{\xi^{(k+1)}}{\lambda^{(k+1)}}\norm{\bx^{(k+1)}-\bx^*}_2.
\end{align}
Combining the bounds for both cases, \eqref{eq:x-bound-alpha} and
\eqref{eq:x-bound-xi}, and using triangular inequality, we have
\begin{align}
\label{eq:x-bound}
\left(\bar{\tau}-\frac{\xi^{(k+1)}}{\lambda^{(k+1)}}\right)\norm{\bx^{(k+1)}-\bx^*}_2\leq
\bar{F}^*+\bar{\tau}\norm{\bx^*}_2 + \frac{\lambda^{(k+1)}}{2}
\left(\norm{\theta^{(k+1)}}_2^2 +
  \frac{\alpha^{(k+1)}}{\left(\lambda^{(k+1)}\right)^2}\right),
\end{align}
for all $k\geq 0$. Note that $\{\lk,\ak,\xik\}$ is chosen in \fal~such that
$\frac{\ak}{(\lk)^2}=\frac{\alpha^{(1)}}{(\lambda^{(1)})^2}$ for all
$k>1$, and both $\frac{\xik}{\lk}\searrow 0$ and $\lk\searrow 0$
monotonically. Since $\sigma_{\min}(A_i)\geq 1$ for all $i\in\cN$, the inductive assumption $\norm{\xk-\bx^*}_2\leq
R$, \eqref{gamma1}, and Lemma~\ref{lem:theta-bound} together imply that
\begin{align}\label{theta-bound}
\norm{\theta^{(k+1)}}_2 \leq
  \min_{i\in\cN} \left\{\max\left\{\sqrt{2 L_i^{(1)}
      \frac{\alpha^{(1)}}{(\lambda^{(1)})^2}},~\frac{\xi^{(1)}}{\lambda^{(1)}}\right\}
  + B_i+\norm{\grad\gamma_i(x^*_i)}_2+L_{\gamma_i}R\right\}.
\end{align}
To simplify bounds further, 
choose $\ai=\tfrac{1}{4N}\left(\li\bar{\tau}\right)^2$, and
$\xii=\tfrac{1}{2}\li\bar{\tau}$ for $\li\leq\sigma^2_{\max}(A)/\bar{L}$, where $\bar{L}=\max_{i\in\cN}\{L_{\gamma_i}\}$. Let $\bar{B}:=\max_{i\in\cN}B_i$ and $\bar{G}:=\max\{\norm{\grad\gamma_i(x_i^*)}_2:\ i\in\cN\}$. Together with
\eqref{eq:x-bound}, \eqref{theta-bound} and $\sigma_{\max}(A)\geq
  1$, this choice of parameters implies that
\begin{align*}
\frac{\bar{\tau}}{2}\norm{\bx^{(k+1)}-\bx^*}_2\leq
\bar{F}^*+\bar{\tau}\norm{\bx^*}_2+ \frac{\lambda^{(1)}}{2} \left[\left(
    \frac{\bar{\tau}\sigma_{\max}(A)}{\sqrt{N}} +
     \bar{B}+\bar{G} + \bar{L}R
  \right)^2 + \frac{\bar{\tau}^2}{4N} \right].
\end{align*}
Define $\beta_1:=\frac{2}{\bar{\tau}}\left(\bar{F}^*+\bar{\tau}\norm{\bx^*}_2\right)$,
$\beta_2:=\frac{\bar{\tau}\sigma_{\max}(A)/\sqrt{N}+\bar{B}+\bar{G}}
{\sqrt{\bar{\tau}}}$,
$\beta_3:=\frac{\bar{L}}{\sqrt{\bar{\tau}}}$, and
$\beta_4:=\frac{\bar{\tau}}{4N}$. Then we have that
$\norm{\bx^{(k+1)} - \bx^*}_2\leq \beta_1 + \lambda^{(1)} \left[
  \Big(\beta_2 + \beta_3 R \Big)^2 + \beta_4 \right]$. 

Note that we are free to choose any $\li>0$ satisfying $\li\leq\sigma^2_{\max}(A)/\bar{L}$. Our objective is to
show that by appropriately choosing $\lambda^{(1)}$, we can guarantee that
$\beta_1 + \lambda^{(1)} \left[  \big(\beta_2 + \beta_3 R \big)^2 +
  \beta_4 \right] \leq R$, which would then complete the inductive
proof. This is indeed true if the above quadratic inequality in $R$,
has a solution, or equivalently if the discriminant $$\Delta = (
2\lambda^{(1)} \beta_2 \beta_3 -1 )^2 - 4\lambda^{(1)} \beta_3^2 \big[
\lambda^{(1)}( \beta_2^2 + \beta_4) + \beta_1\big]$$ is
non-negative. Note that $\Delta$ is continuous in $\lambda^{(1)}$,
and $\lim_{\lambda^{(1)} \rightarrow 0} \Delta = 1$. Thus, for all
sufficiently small $\lambda^{(1)}>0$, we have $\Delta \geq 0$. Hence, we
can set $R=\frac{1 - 2\lambda^{(1)} \beta_2 \beta_3 - \sqrt{\Delta}}{2
  \lambda^{(1)} {\beta_3}^2}$ for some $\lambda^{(1)}>0$ such that
$\Delta\geq 0$, and this will imply that $\norm{\bx^{(k+1)}-\bx^*}_2\leq
R$ whenever $\norm{\xk-\bx^*}_2\leq R$ for all $k\geq1$.

The induction will be complete if we can show that
$\norm{\bx^{(1)}-\bx^*}_2\leq R$. Note that in \fal~we set
$\theta^{(1)}=\mb{0}$. Hence, for $k=0$, \eqref{eq:x-bound} implies that
$\norm{\bx^{(1)}-\bx^*}_2\leq\beta_1+\li\beta_4$. Hence, our choice of $R$
guarantees that $\norm{\bx^{(1)}-\bx^*}_2\leq R$. This completes the
induction.

Following the same arguments leading to \eqref{eq:x-bound}, it can also be shown that for all $k\geq 0$
\begin{equation*}
\left(\bar{\tau}-\frac{\xi^{(k+1)}}{\lambda^{(k+1)}}\right)\norm{\bx_*^{(k+1)}-\bx^*}_2\leq
\bar{F}^*+\bar{\tau}\norm{\bx^*}_2 + \frac{\lambda^{(k+1)}}{2}
\norm{\theta^{(k+1)}}_2^2.
\end{equation*}
Therefore, we can conclude that $\norm{\xkopt-\bx^*}\leq R$ for all $k\geq 1$ holds for the same $R$ we selected above.

Note that $\Delta$ is a concave quadratic of $\lambda^{(1)}$ such that
$\Delta=1$ when $\lambda^{(1)}=0$; hence, one of its roots is positive and
the other one is negative. Moreover, $R\leq\frac{1}{2 \lambda^{(1)}
  {\beta_3}^2} - \frac{\beta_2} {\beta_3}$ and the bound on $R$ is
decreasing in $\lambda^{(1)}>0$. Hence, in order to get a smaller bound on
$R$, we will choose $\lambda^{(1)}$ as the positive root of $\Delta$. In
particular, we set
$\lambda^{(1)}=\frac{\sqrt{\left(\beta_2+\beta_3\beta_1\right)^2+\beta_4}
  -\left(\beta_2+\beta_3\beta_1\right)}{2\beta_3\beta_4}$.
\endproof
\subsection{Proof of Theorem~\ref{thm:iter-complexity}}
\proof
The proof directly follows from Theorem~3.3 in~\cite{ia10}. For the sake of completeness, we also provide the proof here. Let $\bx^*$ denote an optimal solution to \eqref{eq:general_dist}.

Note that (a) follows immediately from Cauchy-Schwarz and the definition of $\theta^{(k+1)}$.
\vspace{-3mm}
\begin{equation*}
\norm{A\xk-b}_2\leq\norm{A\xk-b-\lk\tk}_2+\lk\norm{\tk}_2=\lk(\norm{\theta^{(k+1)}}_2+\norm{\tk}_2)\leq 2B_\theta\lk.
\vspace{-3mm}
\end{equation*}
First, we prove the second inequality in (b). Suppose that $\xk$ satisfies \eqref{eq:inexact}(a), which implies that $\bar{F}(\xk)+\tfrac{\lk}{2}\norm{\theta^{(k+1)}}_2^2\leq\bar{F}(\bx^*)+\tfrac{\lk}{2}\norm{\tk}_2^2+\tfrac{\ak}{\lk}$. Now, suppose that $\xk$ satisfies \eqref{eq:inexact}(b). From the convexity of $\pk$ and Cauchy-Schwarz, it follows that $\pk(\xk)\leq\pk(\bx^*)+\xik\norm{\xk-\bx^*}_2$. Hence, dividing it by $\lk$, we have $\bar{F}(\xk)+\tfrac{\lk}{2}\norm{\theta^{(k+1)}}_2^2\leq\bar{F}(\bx^*)+\tfrac{\lk}{2}\norm{\tk}_2^2+\tfrac{\xik}{\lk}$. Therefore, for all $k\geq 1$, $\xk$ satisfies the second inequality in (b) since it also satisfies
\vspace{-3mm}
\begin{equation*}
\bar{F}(\xk)-\bar{F}^*\leq\lk\left(\frac{\norm{\tk}_2^2-\norm{\theta^{(k+1)}}_2^2}{2}+\frac{\max\left\{\ak ,\xik\norm{\xk-\bx^*}_2\right\}}{(\lk)^2}\right).
\vspace{-3mm}
\end{equation*}
Now, in order to prove the first inequality in (b), we will exploit the primal-dual relations of the following two pairs of problems:
\begin{equation*}
\hspace{-2mm}
\begin{array}{ll}
(\cP): \min_{\bx\in\reals^{nN}}\{\bar{F}(\bx):\ A\bx=b\}, &(\cD):\max_{\theta\in\reals^m}b^\T\theta-\bar{F}^*(A^\T\theta),\\
(\cP_k): \min_{\bx\in\reals^{nN}}\lk\bar{F}(\bx)+\tfrac{1}{2}\norm{A\bx-b_k}_2^2, &(\cD_k):\max_{\theta\in\reals^m}\lk(b^\T\theta-\bar{F}^*(A^\T\theta))-\frac{(\lk)^2}{2}h(\theta),
\end{array}
\end{equation*}
where $b_k:=b+\lk\tk$, $h(\theta):=\norm{\theta-\tk}_2^2-\norm{\tk}_2^2$, and $\bar{F}^*$ denotes the convex conjugate of $\bar{F}$. Note that problem $(\cP_k)$ is nothing but the subproblem in \eqref{eq:subproblem}. Therefore, from weak-duality between $(\cP_k)$ and $(\cD_k)$, it follows that
\vspace{-3mm}
\begin{equation*}
\pk(\xk)=\lk\bar{F}(\xk)+\tfrac{1}{2}\norm{A\xk-b_k}_2^2\geq\lk(b^\T\theta^*-\bar{F}^*(A^\T\theta^*))-\frac{(\lk)^2}{2}h(\theta^*).
\vspace{-3mm}
\end{equation*}
Note that from strong duality between $(\cP)$ and $(\cD)$, it follows that $\bar{F}^*=\bar{F}(\bx^*)=b^\T\theta^*-\bar{F}^*(A^\T\theta^*)$. Therefore, dividing the above inequality by $\lk$, we obtain
\vspace{-3mm}
\begin{equation*}
\bar{F}(\xk)-\bar{F}^*\geq-\frac{\lk}{2}\left(\norm{\theta^*}_2^2-2(\theta^*)^\T\tk+\norm{\theta^{(k+1)}}_2^2\right)\geq-\frac{\lk}{2}\left(\norm{\theta^*}_2+B_\theta\right)^2.
\vspace{-3mm}
\end{equation*}
\endproof
\subsection{Proof of Theorem~\ref{thm:rate}}
\proof
We assume that $\sigma_{\max}(A)\geq \sqrt{\max_{i\in\cN}d_i+1}$, and $\sigma_{\min}(A_i)=\sqrt{d_i}\geq 1$ for all $i\in\cN$, where $d_i$ denotes the degree of $i\in\cN$. As discussed in the proof of
Theorem~\ref{thm:bounded-x}, this is a valid assumption for distributed
optimization problem in \eqref{eq:dist_opt}.
Let $\theta^*$ denote an optimal dual
solution to \eqref{eq:general_dist}. Note that from the first-order
optimality conditions for \eqref{eq:general_dist}, we have
$\mathbf{0}\in\grad\gamma_i(x^*_i)+A_i^\T\theta^*+\partial\rho_i(x_i)|_{x_i=x^*_i}$;
hence, $\norm{A_i^\T\theta^*}_2\leq B_i+ G_i$. 
Therefore, $\norm{\theta^*}_2\leq\min_{i\in\cN}\frac{B_i+G_i}{\sigma_{\min}(A_i)}$.

Given $0<\li\leq\sigma^2_{\max}(A)/\bar{L}$, choose $\ai,\xii>0$ such that $\ai=\tfrac{1}{4N}\left(\li\bar{\tau}\right)^2$, and
$\xii=\tfrac{1}{2}\li\bar{\tau}$.
Then Lemma~\ref{lem:theta-bound} and $\sigma_{\max}(A)\geq 1$ together imply that for
all $k\geq 1$
\begin{equation}
\norm{\tk}_2\leq\min_{i\in\cN}\left\{
\frac{\bar{\tau}\sigma_{\max}(A)/\sqrt{N}+B_i+G_i}{\sigma_{\min}(A_i)}\right\}:=B_\theta.
\end{equation}
Hence, note that $\norm{\theta^*}_2\leq B_\theta$.

To simplify notation, suppose that $\li=\min\left\{1,\sigma^2_{\max}(A)/\bar{L}\right\}=1$. \eqref{eq:x-bound} implies that for all $k\geq 1$
\begin{equation}
\label{eq:Bx}
\norm{\xk-\bx^*}_2 \leq\frac{2}{\bar{\tau}}\left[\bar{F}^*+\bar{\tau}\norm{x^*}_2
  + \tfrac{1}{2}\left(B_\theta^2+\frac{\bar{\tau}^2}{4N}\right)\right]:=B_x.
\end{equation}
Note that \eqref{eq:Bx} implies that
$\frac{\xii}{(\li)^2}B_x=\frac{1}{\li}\frac{\bar{\tau}}{2}B_x\geq\frac{1}{2}
B^2_\theta+\frac{\bar{\tau}^2}{8N}\geq\frac{5}{8N}\bar{\tau}^2\geq\frac{\ai}{(\li)^2}$,
where we used the fact $B_\theta\geq\frac{\sigma_{\max}(A)}{\max_{i\in\cN}\{\sigma_{\min}(A_i)\}}~\frac{\bar{\tau}}{\sqrt{N}}\geq\frac{\bar{\tau}}{\sqrt{N}}$. Note that the last inequality follows from our assumption on $A$ stated at the beginning of the proof, i.e. $\sigma_{\max}(A)\geq\sqrt{\max_{i\in\cN}d_i+1}$ and $\sigma_{\min}(A_i)=d_i$ for all $i\in\cN$. Hence,
Theorem~\ref{thm:iter-complexity}, $\li=1$, and $\norm{\theta^*}_2\leq B_\theta$ imply
that
\begin{align}
N^{f}_{\fal}(\epsilon)&\leq
\log_{\tfrac{1}{c}}\left(\frac{2B_\theta}{\epsilon}\right)
=\log_{\tfrac{1}{c}}\left(2\min_{i\in\cN}\left\{
\frac{\bar{\tau}\sigma_{\max}(A)/\sqrt{N}+B_i+G_i}{\sigma_{\min}(A_i)\epsilon}\right\}\right):=\bar{N}^f,\label{eq:fal_iter_feasibility}\\
N^{o}_{\fal}(\epsilon)&\leq
\log_{\tfrac{1}{c}}\left(\frac{1}{\epsilon}~\max
  \left\{\tfrac{1}{2}\left(\norm{\theta^*}_2+B_\theta\right)^2,
~B_\theta^2+\bar{F}^*+\bar{\tau}\norm{x^*}_2+\frac{\bar{\tau}^2}{8N}\right\}\right),\nonumber\\
&=\log_{\tfrac{1}{c}}\left(\frac{2B_\theta^2+\bar{F}^*+\bar{\tau}\norm{x^*}_2+\frac{\bar{\tau}^2}{8N}}{\epsilon}\right)
:=\bar{N}^o\label{eq:fal_iter_optimality}.
\end{align}

Since $\ai=\tfrac{1}{4N}\left(\li\bar{\tau}\right)^2$, we have
$\sqrt{\ak}=\frac{\bar{\tau}}{\sqrt{4N}} c^k$. Hence, Lemma~\ref{lem:subiter-complexity} 
implies that
\begin{equation}
\label{eq:apg_iter}
\Nk\leq 2B_x\sqrt{\frac{2(\lk\bar{L}+\sigma^2_{\max}(A))}{\ak}}\leq \frac{8B_x
  \sqrt{N}}{\bar{\tau}}\sigma_{\max}(A) c^{-k}.
\end{equation}
Hence, \eqref{eq:fal_iter_feasibility} and \eqref{eq:apg_iter} imply that
the total number of \apg~iterations to compute an $\epsilon$-feasible
solution can be bounded above:
\begin{align}
\hspace{-2mm}
\sum_{k=1}^{N^{f}_{\fal}(\epsilon)}\Nk&\leq \frac{8B_x
  \sqrt{N}}{\bar{\tau}}\sigma_{\max}(A)\sum_{k=1}^{\bar{N}^f}c^{-k}\leq
\frac{8B_x \sqrt{N}}{c(1-c)\bar{\tau}}\sigma_{\max}(A)
\left(\frac{1}{c}\right)^{\bar{N}^f},\nonumber\\
&\leq \frac{16B_x
  \sqrt{N}}{c(1-c)\bar{\tau}}\min_{i\in\cN}\left\{
\frac{\bar{\tau}\sigma_{\max}(A)/\sqrt{N}+B_i+G_i}{\sigma_{\min}(A_i)\epsilon}\right\}\frac{\sigma_{\max}(A)}{\epsilon}=\cO\left(\frac{\sigma^2_{\max}(A)}{\min_{i\in\cN}\sigma_{\min}(A_i)}\frac{1}{\epsilon}\right). \nonumber
\end{align}

Similarly, \eqref{eq:fal_iter_optimality} and \eqref{eq:apg_iter} imply
that the total number of \apg~iterations to compute an $\epsilon$-optimal
solution can be bounded above:
\begin{align}
\sum_{k=1}^{N^{o}_{\fal}(\epsilon)}\Nk\leq \frac{8B_x
  \sqrt{N}}{c(1-c)\bar{\tau}}\sigma_{\max}(A)
\left(\frac{1}{c}\right)^{\bar{N}^0}=\cO\left(\frac{\sigma^3_{\max}(A)}{\min_{i\in\cN}\sigma^2_{\min}(A_i)}\frac{1}{\epsilon}\right).
\end{align}
\endproof
\subsection{Proof of Lemma~\ref{lem:subgrad}}
\proof
Given any convex function $\rho:\reals^n\rightarrow\reals$ and $\bar{x}\in\reals^n$, in order to simplify the notation throughout the proof, $\partial\rho(x)|_{x=\bar{x}}\subset\reals^n$, the subdifferential of $\rho$ at $\bar{x}$, will be written as $\partial\rho(\bar{x})$. Given $\bar{x}\in\reals^n$, there exists $\nu\in\partial P(\bar{x})$ such that $\norm{\nu}_2\leq\xi$, if and only if $\norm{\nu^*}\leq\xi$, where $\nu^*=\argmin\{\norm{\nu}_2:\ \nu\in\partial P(\bar{x})\}$. Note that $\partial P(\bar{x})=\lambda \partial\rho(\bar{x})+\grad f(\bar{x})$, and
\begin{equation}
\partial\rho(\bar{x})= \beta_1\prod_{k=1}^K \partial \norm{\bar{x}_{g(k)}}_1 +\beta_2 \prod_{k=1}^K \partial \norm{\bar{x}_{g(k)}}_2,
\end{equation}
where $\prod$ denotes the Cartesian product. Since the groups $\{g(k)\}_{k=1}^K$ are not overlapping with each other, the minimization problem is separable in groups. Hence, for all $k\in [1,K]$, we have $\nu^*_{g(k)}=\pi^*_{g(k)}+\omega^*_{g(k)}+\grad_{x_{g(k)}}f(\bar{x})$ such that
\begin{equation}
\label{eq:subgrad_check}
\begin{array}{ll}
(\pi^*_{g(k)},\omega^*_{g(k)})&=\argmin \norm{\pi_{g(k)}+\omega_{g(k)}+\grad_{x_{g(k)}}f(\bar{x})}_2^2\\ &s.t.\quad \pi_{g(k)}\in\lambda\beta_1\partial \norm{\bar{x}_{g(k)}}_1, \quad \omega_{g(k)}\in\lambda\beta_2\partial \norm{\bar{x}_{g(k)}}_2.
\end{array}
\end{equation}
Fix $k\in[1,K]$. We will consider the solution to above problem in two cases. Suppose that $\bar{x}_{g(k)}=\mathbf{0}$. Since $\partial\norm{\mathbf{0}}_1$ is the unit $\ell_\infty$-ball, and $\partial\norm{\mathbf{0}}_2$ is the unit $\ell_2$-ball, \eqref{eq:subgrad_check} can be equivalently written as
\begin{equation}
\label{eq:subgrad_check-0}
\begin{array}{ll}
(\pi^*_{g(k)},\omega^*_{g(k)})&=\argmin \norm{\pi_{g(k)}+\omega_{g(k)}+\grad_{x_{g(k)}}f(\bar{x})}_2^2\\
&s.t.\quad \norm{\pi_{g(k)}}_\infty\leq\lambda\beta_1, \quad \norm{\omega_{g(k)}}_2\leq\lambda\beta_2.
\end{array}
\end{equation}
Clearly, it follows from Euclidean projection on to $\ell_2$-ball that $$\omega^*_{g(k)}=-(\pi^*_{g(k)}+\grad_{x_{g(k)}}f(\bar{x}))
\min\left\{1,\frac{\lambda\beta_2}{\norm{\pi^*_{g(k)}+\grad_{x_{g(k)}}f(\bar{x})}_2}\right\}.$$
Hence, $\norm{\pi^*_{g(k)}+\omega^*_{g(k)}+\grad_{x_{g(k)}}f(\bar{x})}_2=\max\{0,~\norm{\pi^*_{g(k)}+\grad_{x_{g(k)}}f(\bar{x})}_2-\lambda\beta_2\}$. Therefore,
\begin{align*}
\pi^*_{g(k)}=\argmin\{\norm{\pi_{g(k)}+\grad_{x_{g(k)}}f(\bar{x})}_2:\ \norm{\pi_{g(k)}}_\infty\leq\lambda\beta_1\}=-\sgn(\grad_{x_{g(k)}}f(\bar{x}))\odot\min\{|\grad_{x_{g(k)}}f(\bar{x})|, \lambda\beta_1\}.
\end{align*}

Now, suppose that $\bar{x}_{g(k)}\neq\mathbf{0}$. This implies that $\partial \norm{\bar{x}_{g(k)}}_2=\{\bar{x}_{g(k)}/\norm{\bar{x}_{g(k)}}_2\}$. Hence, when $\bar{x}_{g(k)}\neq\mathbf{0}$, we have $\omega^*_{g(k)}=\lambda\beta_2\bar{x}_{g(k)}/\norm{\bar{x}_{g(k)}}_2$, and the structure of $\partial\norm{\cdot}_1$ implies that $\pi^*_j=\lambda\beta_1\sgn\left(\bar{x}_j\right)$ for all $j\in g(k)$ such that $|\bar{x}_j|>0$; and it follows from \eqref{eq:subgrad_check} that for all $j\in g(k)$ such that $\bar{x}_j=0$, we have
\begin{align*}
\pi^*_j=\argmin\left\{\left(\pi_j+\tfrac{\partial}{\partial x_j}f(\bar{x})\right)^2:\ |\pi_j|\leq\lambda\beta_1\right\}=-\sgn\left(\tfrac{\partial}{\partial x_j}f(\bar{x})\right)
\min\left\{\left|\tfrac{\partial}{\partial x_j}f(\bar{x})\right|,\ \lambda\beta_1\right\}.
\end{align*}
\endproof
\vspace{-1cm}
\subsection{Proof of Lemma~\ref{lem:prox}}
\proof
Since the groups are not overlapping with each other, the proximal problem
becomes separable in groups. Let $n_k:=|g(k)|$ for all $k$. Thus, it suffices to show that $\min_{x_{g(k)}
  \in \reals^{n_k}} \{ \beta_1 \norm{x}_1 + \beta_2 \norm{x_{g(k)}}_2 +
\frac{1}{2t} \norm{ x_{g(k)} - \bar{x}_{g(k)}}_2^2  \}$ has a closed form
solution as shown in the statement for some fixed $k$. By the definition
of dual norm, we have
\begin{align}
&\min_{x_{g(k)} \in \reals^{n_k}} \beta_1 \norm{x_{g(k)}}_1 + \beta_2 \|
x_{g(k)} \|_2 + \frac{1}{2t} \| x_{g(k)} - \bar{x}_{g(k)}
\|_2^2, \label{eq:prox_x}\\
= & \min_{x_{g(k)} \in \mathbb{R}^{n_k}} \: \max_{ \| u_1\|_{\infty} \leq
  \beta_1} u_1^{\mathsf{T}} x_{g(k)} + \max_{ \| u_2\|_{2} \leq \beta_2
  } u_2^{\mathsf{T}} x_{g(k)} + \frac{1}{2t} \| x_{g(k)} - \bar{x}_{g(k)}
\|_2^2, \nonumber\\
= & \max_{\substack{
    \| u_1\|_{\infty} \leq \beta_1\\
    \| u_2\|_{2} \leq \beta_2 }}
\min_{x \in \mathbb{R}^n} (u_1 + u_2)^{\mathsf{T}} x_{g(k)} + \frac{1}{2t}
\| x_{g(k)} - \bar{x}_{g(k)} \|_2^2, \label{eq:prox_saddle}\\
= & \max_{\substack{
    \| u_1\|_{\infty} \leq \beta_1\\
    \| u_2\|_{2} \leq \beta_2 }}
(u_1 + u_2)^{\mathsf{T}} \bar{x}_{g(k)} - \frac{t}{2} \| u_1 + u_2
\|_2^2. \label{eq:prox_u}
\end{align}
Let $(u_1^*,u_2^*)$ be the optimal solution of \eqref{eq:prox_u}. Since
$x^p_{g(k)}$ is the optimal solution to \eqref{eq:prox_x}, it follows from
\eqref{eq:prox_saddle} that
\begin{align}\label{x*}
x^p_{g(k)} = \bar{x}_{g(k)} - t(u^*_1 + u^*_2).
\end{align}
Note that \eqref{eq:prox_u} can be equivalently written as $\min\{\| u_1 +
u_2 - \frac{1}{t} \bar{x}_{g(k)} \|_2^2:\ \| u_1\|_{\infty} \leq \beta_1,
\| u_2\|_{2} \leq \beta_2\}$. Minimizing over $u_2$, we have
\begin{align}\label{u2}
u^*_2(u_1) = \left( \frac{1}{t}\bar{x}_{g(k)} -  u_1\right)\min \left\{
\frac{\beta_2} { \| \frac{1}{t}\bar{x}_{g(k)} - u_1 \|_2} , 1
\right\}.
\end{align}
Hence, we have \vspace{-5mm}
\begin{align*}
u_1^*&=\argmin_{\| u_1\|_{\infty} \leq \beta_1} \left \| \left( u_1
    -\frac{1}{t}\bar{x}_{g(k)}\right) \max \middle\{  1- \frac{\beta_2 }
  { \|u_1 - \frac{1}{t}\bar{x}_{g(k)}\|_2}, 0  \middle\} \right\|_2=\argmin_{\| u_1\|_{\infty} \leq \beta_1}\max\{\norm{u_1 -
  \frac{1}{t}\bar{x}_{g(k)}}_2-\beta_2 ,~0\}.
\end{align*}
Clearly, $u_1^*=\argmin_{\| u_1\|_{\infty} \leq \beta_1}  \| ( u_1
-\frac{1}{t}\bar{x}_{g(k)}) \|_2=\sgn(\bar{x}_{g(k)}) \min \left\{
  \frac{1}{t}| \bar{x}_{g(k)} |, \beta_1  \right\}$.
The final result follows from combining \eqref{x*} and \eqref{u2}.
\endproof
\subsection{Improved rate for asynchronous DFAL}
\label{sec:abcd}
Let $\cR$ denote a discrete  random variable uniformly distributed over
the set $\cN$. 
Let $[U_1,U_2,\ldots,U_N]$ denote a partition of the
$nN$-dimensional identity matrix where $U_i \in \reals^{nN\times n}$, $i =
1, \ldots, N$.  
In the rest, given $\bh\in\reals^{nN}$, we denote
$\bh_{[\cR]}:=U_\cR U_{\cR}^\top \bh$. Consider the composite convex optimization
problem
\begin{equation}
\label{eq:composite_problem}
\Phi^*:=\min_{\by \in \mathbb{R}^{nN}} \Phi(\by) := \sum_{i=1}^N \rho_i(y_i) +f(\by),
\end{equation}
where $\rho_i:\reals^n\rightarrow \reals$ is a closed convex function for
all $i\in\cN$ such that $\prox{t\rho_i}$ can be computed efficiently for
all $t>0$ and $i\in\cN$, and $f:\reals^{nN}\rightarrow\reals$ is a
differentiable convex function such that for some
$\{L_i\}_{i\in\cN}\subset\reals_{++}$, $f$ satisfies
\begin{equation}
\mathbb{E}[f(\by+\bh_{[\cR]})]\leq f(\by)+\frac{1}{N}\left(\fprod{\grad
    f(\by),\bh}+\frac{1}{2}\sum_{i\in\cN}L_i\norm{h_i}_2^2\right) \label{eq:ESO}
\end{equation}
for all $\by,\bh\in\reals^{nN}$.
\citet{fercoq2013accelerated} proposed the accelerated proximal coordinate descent
algorithm \abcd\ (see Figure~\ref{alg:abcd}) to solve
\eqref{eq:composite_problem}. They showed that for a given
$\alpha>0$, the iterate sequence $\{\bz^{(\ell)},\bu^{(\ell)}\}$
computed by \abcd\ 
satisfies \vspace{-0.15cm}
\begin{equation}
\mathbb{E}\left[
  \Phi\left(\left(\frac{1}{Nt^{(\ell)}}\right)^2\bu^{(\ell+1)}
    +\bz^{(\ell+1)}\right)
  - \Phi^*\right] \leq \alpha,\quad \forall \ell\geq
2N\sqrt{\frac{C}{\alpha}}, \label{eq:abcd_complexity}
\end{equation}
where
\begin{equation}
  C:=\min_{\by^*\in\cY^*}(1-\tfrac{1}{N})\left(\Phi\left(\bz^{(0)}\right)-
    \Phi^*\right)
  +\tfrac{1}{2}\sum_{i\in\cN}L_i\norm{z_i^{(0)}-y^*_i}_2^2, \label{eq:C}
\end{equation}
and $\cY^*$ denotes the set of optimal solutions.

\begin{figure}[h!]
    \rule[0in]{6.2in}{1pt}\\
    \textbf{Algorithm ARBCD}$~(\bz^{(0)})$\\
    \rule[0.125in]{6.2in}{0.1mm}
    \vspace{-0.75cm}
    {\footnotesize
    \begin{algorithmic}[1]
    \STATE $\ell\gets0,\quad t^{(0)}\gets 1,\quad u_i^{(1)}\gets
    \mathbf{0},\quad \forall i\in\cN$
    \WHILE{$\ell\geq 0$}
    \STATE 
    $i$ is a sample of $\cR$
    \STATE
    $z_i^{(\ell+1)}\gets\prox{t^{(\ell)}\rho_i/L_i}\left(z_i^{(\ell)}-\frac{t^{(\ell)}}{L_i}\grad_{y_i}f\left(\left(\frac{1}{N
            t^{(\ell)}}\right)^2\bu^{(\ell)}+\bz^{(\ell)}\right)\right)$
    \STATE $u_i^{(\ell+1)}\gets u_i^{(\ell)}+N^2 t^{(\ell)}(1-t^{(\ell)})\left(z_i^{(\ell+1)}-z_i^{(\ell)}\right)$
    \STATE $z_{-i}^{(\ell+1)}\gets z_{-i}^{(\ell)},\quad u_{-i}^{(\ell+1)}\gets u_{-i}^{(\ell)}$
    \STATE $t^{(\ell+1)}\gets \frac{1+\sqrt{1+(2Nt^{(\ell)})^2}}{2N}$
    \ENDWHILE
    \end{algorithmic}
    \rule[0.25in]{6.2in}{0.1mm}
    }
    \vspace{-0.75cm}
    \caption{\small Accelerated Randomized Proximal Block Coordinate Descent~(ARBCD) algorithm}
    \label{alg:abcd}
\end{figure}
In the following result, we establish that the bound~\eqref{eq:ESO}
can be exploited for designing an accelerated version of
asynchronous \fal.
\begin{lemma}
\label{lem:abcd}
 Fix $\alpha>0$, and $p\in(0,1)$. Let
 $\{\bz_k^{(\ell)},\bu_k^{(\ell)}\}_{\ell\in\integers_+}$, $k = 1, \ldots,
 K$, denote the iterate sequence corresponding to $K:=\log(1/p)$
 independent calls to \abcd$(\by^{(0)})$.
 Define
 $\by_k:=\left(\frac{1}{Nt^{(T)}}\right)^2\bu_k^{(T+1)}+\bz_k^{(T+1)}$ for
 $k = 1, \ldots, K$, and $T:=2N\sqrt{\frac{2\cC}{\alpha}}$.
 Then
 \eq
 \mathbb{P}\left(\min_{k=1,\ldots,K}\Phi(\by_k)-\Phi^*\leq\alpha\right)\geq
 1-p.
 \en
\end{lemma}
\proof Since the sequence $\{\by_k\}_{k=1}^K$ is i.i.d., and each
$\by_k$ satisfies $\mathbb{E}[\Phi(\by_k)-\Phi^*]\leq\tfrac{\alpha}{2}$,
Markov's inequality implies that
$\mathbb{P}(\Phi(\by_k)-\Phi^*>\alpha)\leq\mathbf{E}[\Phi(\by_k)
-\Phi^*]/\alpha\leq
\frac{1}{2}$ for $1\leq k\leq K$. Therefore, we have
\begin{equation*}
\mathbb{P}\left(\min_{k=1,\ldots,K}\Phi(\by_k)-\Phi^*\leq\alpha\right)=1-\prod_{k=1}^K\mathbb{P}(\Phi(\by_k)-\Phi^*>\alpha)\leq \left(\frac{1}{2}\right)^K=1-p.
\end{equation*}
\qed

From Lemma~\ref{lem:abcd} 
it follows that we can compute $\by_\alpha$ such that
$\mathbb{P}\left(\Phi(\by_\alpha)-\Phi^*\leq\alpha\right)\geq 1-p$ in at
most
$2N\sqrt{\frac{2C}{\alpha}}\log(\frac{1}{p})$  \abcd~iterations. 
This new oracle can be used to construct an asynchronous version of
\fal~algorithm  with $\cO(1/\epsilon)$
complexity. 
\begin{thm}
Fix $\epsilon > 0 $ and $p\in(0,1)$.
Consider a \emph{asynchronous} variant of \fal\ 
where
\eqref{eq:inexact}(a) in
Figure~\ref{fig:fal} is replaced by
\begin{equation}
\label{eq:inexact_asynch_2}
\mathbb{P} \left( P^{(k)} \big( \mathbf{x}^{(k)} \big) -  P^{(k)} \big(
  \mathbf{x}_*^{(k)} \big) \leq \ak \right) \geq \big(1 - p
\big)^{\frac{1}{N(\epsilon)}},
\end{equation}
where $N(\epsilon) =
\log_{\frac{1}{c}}\left(\frac{\bar{C}}{\epsilon}\right)$ is defined in
Corollary~\ref{cor:N_eps}. Then 
$\{x_i^{\left(N(\epsilon)\right)}\}_{i \in \cN}$, satisfies
\begin{equation*}
  \mathbf{P}(\epsilon):=\mathbb{P} \left(
    \big|\sum_{i\in\cN}F_i\left(x^{\left(N(\epsilon)\right)}_i\right)
    -F^*\big|\leq\epsilon,\
    \hbox{ and }
    \max_{(i,j)\in\cE}\big\{\norm{x^{\left(N(\epsilon)\right)}_i
      -x^{\left(N(\epsilon)\right)}_j}_2\big\}\leq\epsilon\right)
  \geq 1-p,
  \vspace{-0.25cm}
\end{equation*}
and $\cO\left(\frac{1}{\epsilon}\log\left(\frac{1}{p}\right)\right)$
\abcd~iterations are required to
compute $\{x_i^{\left(N(\epsilon)\right)}\}_{i \in \cN}$.
\proof
Consider the $k$-th \fal~subproblem $\min P^{(k)}(\bx):=\lk
\sum_{i\in\cN}\rho_i(x_i)+f^{(k)}(\bx)$, where $f^{(k)}$ is defined in
\eqref{eq:fk}. Let $\tilde{L}^{(k)}_i:=\lk L_{\gamma_i}+d_i$ for all
$i\in\cN$. Then it can be easily shown that $f^{(k)}$ satisfies
\eqref{eq:ESO} with constants $\{\tilde{L}^{(k)}_i\}_{i\in\cN}$ for all $1\leq k\leq
N(\epsilon)$. Hence, \abcd~algorithm can be used to solve $\min
P^{(k)}(\bx)$ with the iteration complexity given in
Lemma~\ref{lem:abcd}. Consider the random event
\begin{equation}
\Delta:=\bigcap_{k=1}^{N(\epsilon)}\left\{\pk(\xk)-\pk(\xkopt)\leq\ak\quad \mbox{\textbf{or}}\quad \exists \gk_i\in\partial_{x_i} \pk(\bx)|_{\bx=\xk}\ \mbox{ s.t. }
  \max_{i\in\cN}\norm{\gk_i}_2\leq\tfrac{\xik}{\sqrt{N}}\right\}.
\end{equation}
Clearly, for all random sequences $\{\bx^{(k)}\}_{k=1}^{N(\epsilon)}$
satisfying random event $\Delta$, Corollary~\ref{cor:N_eps} implies that
$\big|\sum_{i\in\cN}F_i\left(x^{\left(N(\epsilon)\right)}_i\right)
-F^*\big|\leq\epsilon$
and
$\max_{(i,j)\in\cE}\big\{\norm{x^{\left(N(\epsilon)\right)}_i
  -x^{\left(N(\epsilon)\right)}_j}_2\big\}\leq\epsilon$. Hence,
we have \vspace{-0.2cm}
\begin{equation*}
\mathbf{P}(\epsilon)\geq
\mathbb{P}(\Delta)\geq\prod_{k=1}^{N(\epsilon)}\mathbb{P} \left( P^{(k)}
  \big( \mathbf{x}^{(k)} \big) -  P^{(k)} \big(
  \mathbf{x}_*^{(k)} \big) \leq \ak \right)\geq 1-p.
\end{equation*}

\vspace{-0.25cm}
\noindent In the rest, we bound the total number of \abcd~iterations
required by \emph{asynchronous} variant of \fal~to compute
$\bx^{(N(\epsilon))}$. Note that $(1-p)^{\frac{1}{N(\epsilon)}}$ is a
concave function for $p\in(0,1)$, and we have
$(1-p)^{\frac{1}{N(\epsilon)}}\leq 1-\frac{p}{N(\epsilon)}$. Therefore,
Lemma~\ref{lem:abcd} and the discussion after Lemma~\ref{lem:abcd}
together imply that the number of \abcd~iterations, $N^{(k)}$, to compute $\bx^{(k)}$ satisfying either
\eqref{eq:inexact_asynch_2} or \eqref{eq:inexact}(b) is bounded above for $1\leq k \leq
N(\epsilon)$ as follows
\begin{equation}
N^{(k)} \leq 2N\sqrt{\frac{2C^{(k)}}{\ak}}\log
\left(\frac{N(\epsilon)}{p}\right)=2N\left(\log
  \left(\frac{1}{p}\right)+\log\log_{\frac{1}{c}}
  \left(\frac{\bar{C}}{\epsilon}\right)
\right)\sqrt{\frac{2C^{(k)}}{\ak}},
\end{equation}
with $C^{(k)} = P^{(k)}\left(\bx^{(k-1)}\right)-P^{(k)} \left(
  \bx_*^{(k)}\right)
+\sum_{i\in\cN}\frac{\tilde{L}^{(k)}_i}{2}\norm{x_i^{(k-1)}-x_{*i}^{(k)}}_2^2$.

Convexity of $\{\rho_i\}_{i\in\cN}$, and Lemma~\ref{lem:lipschitz-fk}
imply that \vspace{-0.25cm}
\begin{equation*}
P^{(k)}(\bx^{(k-1)})-P^{(k)}(\bx_*^{(k)})\leq \fprod{\lk s^{(k)}+\grad f^{(k)}(\bx_*^{(k)}),\ \bx^{(k-1)}-\bx_*^{(k)}}+\sum_{i\in\cN}\frac{L_i^{(k)}}{2}\norm{x_i^{(k-1)}-x_{*i}^{(k)}}_2^2,
\end{equation*}

\vspace{-0.5cm}
\noindent where $s^{(k)}\in\partial\lk\bar{\rho}(\bx)|_{\bx=\bx^{(k-1)}}$, and $\bar{\rho}(\bx)=\sum_{i\in\cN}\rho_i(x_i)$. Note that optimality conditions imply that $-\grad f^{(k)}(\bx_*^{(k)})\in\partial\lk\bar{\rho}(\bx)|_{\bx=\bx_*^{(k)}}$. Assumption~\ref{assump} implies that $\norm{\grad_{x_i} f^{(k)}(\bx_*^{(k)})}_2\leq \lk B_i$ and $\norm{s_i^{(k)}}_2\leq \lk B_i$ for all $i\in\cN$. Hence, for some $\tilde{C}>0$, we have $C^{(k)}\leq\sum_{i\in\cN}\left(\frac{\Lk_i+\tilde{L}^{(k)}_i}{2}+2\lk B_i\right)\norm{x_i^{(k-1)}-x_{*i}^{(k)}}_2^2\leq \tilde{C}B^2_x$ for all $k\geq 1$. Consequently, we can bound the total number of \abcd~iterations to compute $\bx^{(N(\epsilon))}$ as follows:
\begin{equation*}
\sum_{k=1}^{N(\epsilon)}N^{(k)}\leq2N B_x\sqrt{\frac{2\tilde{C}}{\alpha^{(0)}}}\left(\log\left(\frac{1}{p}\right)+\log\log_{\frac{1}{c}}\left(\frac{\bar{C}}{\epsilon}\right)\right)\sum_{k=1}^{N(\epsilon)}c^{-k}.
\end{equation*}
Since $N(\epsilon)=\log_{\frac{1}{c}}(\bar{C}/\epsilon)$, and $\sum_{k=1}^{N(\epsilon)}c^{-k}=\frac{\left(\frac{1}{c}\right)^{N(\epsilon)}-1}{1-c}=\bar{C}\epsilon^{-1}/(1-c)$. Hence, we can conclude that $\sum_{k=1}^{N(\epsilon)}N^{(k)}=\cO\left(\frac{1}{\epsilon}\left(\log\left(\frac{1}{p}\right)+\log\log\left(\frac{1}{\epsilon}\right)\right)\right)$
%
\end{thm}
\end{document}